\documentclass[11pt, a4paper, english]{amsart}

\usepackage[utf8]{inputenc}
\usepackage[T1]{fontenc}
\usepackage[english]{babel}
\usepackage[shortlabels]{enumitem}
\setlist[enumerate]{font={\bfseries}}
\setlist[enumerate]{label={(\roman*)}}

\usepackage[pdftex]{graphicx} 
\graphicspath{{images/}}
\usepackage{booktabs} 
\usepackage{latexsym}
\usepackage{amsmath,amssymb,amsthm,amsfonts,amscd}
\usepackage{mathtools}
\usepackage{mathrsfs}
\usepackage{dsfont}
\usepackage{comment}
\usepackage[dvipsnames]{xcolor}

\usepackage{thmtools}
\usepackage[arrow, matrix, curve]{xy}

\usepackage{leftidx,tensor}

\usepackage{stmaryrd}
\SetSymbolFont{stmry}{bold}{U}{stmry}{m}{n}

\renewcommand{\epsilon}{\varepsilon}

\usepackage{mleftright}

\usepackage{tikz}	
\usepackage{pgfplots}	
\usetikzlibrary{arrows, angles, quotes, calc,through,backgrounds,matrix,decorations.pathmorphing,
	intersections, pgfplots.fillbetween, patterns}
 \usepgfplotslibrary{polar}
\pgfplotsset{compat=newest}
\usepackage{multirow}
\usepgfplotslibrary{fillbetween}
\pgfplotsset{compat=1.17}
\pgfdeclarelayer{ft}
\pgfdeclarelayer{bg}
\pgfsetlayers{bg,main,ft}
\usepackage{tikz-cd}
\usepackage[font=small, format=plain, labelfont=bf, up, justification=justified,
singlelinecheck=false, bf, hang]{caption}
\usepackage{threeparttable}%

\usepackage{tocbasic}

\usepackage{xfrac}

\usepackage{colonequals}

\usepackage{setspace}
\usepackage{multicol}
\usepackage[linkcolor=black, urlcolor=black ,citecolor=black, colorlinks=true]{hyperref}
\usepackage[capitalize, nameinlink, noabbrev]{cleveref} 

\usepackage[square, sort, numbers]{natbib}

\usepackage{units}

\usepackage{kantlipsum} 
\usepackage{xcolor}
\usepackage{csquotes}
\usepackage{nameref}
\usepackage{breqn}
\usepackage{subcaption}
\usepackage[export]{adjustbox}
\usepackage{imakeidx}
\usepackage{adjustbox}
\usepackage{doi}

\usepackage{bm}

\declaretheorem[style=plain, numberwithin=section]{theorem}
\declaretheorem[sibling=theorem]{lemma}
\declaretheorem[sibling=theorem]{corollary}

\declaretheorem[style=definition, sibling=theorem]{definition}

\declaretheorem[style=definition, sibling=theorem]{remark}
\declaretheorem[style=definition, sibling=theorem]{example}

\declaretheorem[name=Theorem]{introtheorem}

\newcounter{introeq}
\renewcommand{\theintroeq}{\arabic{introeq}}  

\newcommand{\im}{\operatorname{Im}}

\newcommand{\cF}{\mathcal{F}}

\newcommand{\cJ}{\mathcal{J}}

\newcommand{\cM}{\mathcal{M}}

\newcommand{\cP}{\mathcal{P}}

\newcommand{\cS}{\mathcal{S}}

\newcommand{\bC}{\mathbb{C}}

\newcommand{\bE}{\mathbb{E}}

\newcommand{\bN}{\mathbb{N}}

\newcommand{\bP}{\mathbb{P}}

\newcommand{\bR}{\mathbb{R}}

\newcommand{\bT}{\mathbb{T}}

\newcommand{\bZ}{\mathbb{Z}}

\newcommand{\bone}{\mathds{1}}

\newcommand\subsetsim{\mathrel{%
\ooalign{\raise0.2ex\hbox{$\subset$}\cr\hidewidth\raise-0.8ex\hbox{\scalebox{0.9}{$\sim$}}\hidewidth\cr}}}

\newcommand{\norm}[1]{\lVert#1\rVert}
\newcommand{\abs}[1]{\lvert#1\rvert}

\DeclareMathOperator{\Spec}{Spec}

\DeclareMathOperator{\supp}{supp}

\DeclareMathOperator{\linspan}{span}

\DeclareMathOperator{\Vol}{Vol}

\DeclareMathOperator{\Var}{Var}

\DeclareMathOperator{\Cov}{Cov}

\DeclareFontFamily{U}{mathx}{}
\DeclareFontShape{U}{mathx}{m}{n}{<-> mathx10}{}
\DeclareSymbolFont{mathx}{U}{mathx}{m}{n}
\DeclareMathAccent{\widehat}{0}{mathx}{"70}
\DeclareMathAccent{\widecheck}{0}{mathx}{"71}
\newcommand\equidist{\overset{d}{=}}

\numberwithin{equation}{subsection}

\begin{document}

\setlist[itemize]{itemsep=2pt, topsep=3pt}
\setlist[enumerate]{label=(\roman*), topsep=3pt, itemsep=2pt}
	
\newcommand{\Mfull}{\ensuremath{M_{\mathrm{full}}}}
\newcommand{\Msph}{\ensuremath{M_{\mathrm{sph}}}}

\title{Hyperuniformity of self-similar point processes}

\author{Daniel Roca Gonz\'alez}
\address{Institute of Algebra and Geometry \\ Karlsruhe Institute of Technology, Germany}
\email{daniel.roca.gonzalez@proton.me}
\urladdr{https://www.math.kit.edu/iag8/~roca/de} 

\begin{abstract}
We study hyperuniformity of self-similar point processes arising from substitution rules in two dimensions.
In particular, we derive a sufficient condition for hyperuniformity of these point processes
only in terms of the associated substitution matrix.
This condition applies to a wide class of examples for which hyperuniformity had not yet been established,
including most well-known examples of planar self-similar tilings.
In particular, we show that the Godr\`eche--Lan\c{c}on--Billard substitution rule
gives rise to hyperuniform point processes with singular continuous diffraction.
Furthermore, we prove that hyperuniformity is not an MLD invariant,
contradicting standing conjectures.
\end{abstract}

\maketitle
\tableofcontents
\section{Introduction}
\subsection{Hyperuniformity}

In this article,
we study hyperuniformity for \emph{self-similar point processes},
arising from the field of aperiodic order.

A point process is said to be \emph{hyperuniform}
if it exhibits anomalously low density fluctuations at large scales:
more precisely,
a locally square-integrable stationary point process \(\Lambda\) on \(\bR^d\) is \emph{hyperuniform}
if \(\lim_{R \to \infty}\frac{\Var(\# \Lambda \cap B_R)}{R^d} = 0\),
where $B_R$ is the Euclidean ball of radius $R$ centered at the origin.
Hyperuniformity was introduced in 2003 by \citet*{Torquato_2003}:
this property has been widely studied from the viewpoints of materials science,
chemistry and stochastic geometry.
See \citet*{Torquato_2018} and \citet*{coste_2021} for recent surveys of the current state of the art.

While there exists a large body of work on hyperuniformity,
the literature is richest
for structures that lack aperiodic order,
such as physical models of exotic gases, liquids and plasmas.
The only rigorous results on hyperuniformity
for structures with aperiodic order
that we are aware of
are those of \citet*{hartnick2022hyperuniformity},
who studied hyperuniformity for cut-and-project sets,
as well as those of \citet*{Baake_2019_renormalisation},
who proved hyperuniformity for several one-dimensional substitution rules.
There are also numerical results and heuristic arguments in this direction
\citep*{Zachary2009, Lin_2017, Ouz2019, torquato-penrose-numeric, Fuchs2019}.

In particular, \citet*{hartnick2022hyperuniformity}
proved
that there exist cut-and-project quasicrystals which are not hyperuniform.
This is surprising,
as quasicrystals are usually thought to be models of materials
with aperiodic order;
one would expect highly ordered structures to have low density fluctuations,
but it turns out that this is not always the case.
This underscores the importance of having rigorous criteria to prove or disprove hyperuniformity
for structures with aperiodic order.

We will focus on \emph{self-similar point processes},
which arise from substitution rules
as we will explain later in this section.
Extending the methods of \citet*{Baake_2019_renormalisation},
we prove hyperuniformity for a wide class of point processes in one and two dimensions,
including the following.
\begin{introtheorem}\label{theorem:examples-intro}
The point processes associated to the following substitution rules are hyperuniform:
\begin{enumerate}
    \item \href{https://tilings.math.uni-bielefeld.de/substitution/penrose-kite-dart/}{Penrose} \citep*{thesis}.
    \item \href{https://tilings.math.uni-bielefeld.de/substitution/chair/}{Chair} \citep*{thesis}.
    \item \href{https://tilings.math.uni-bielefeld.de/substitution/danzers-7-fold/}{Danzer's 7-fold} \citep*{encyclopedia_danzers_7fold}.
    \item \href{https://tilings.math.uni-bielefeld.de/substitution/binary}{Godr\`eche--Lan\c{c}on--Billard} \citep*{encyclopedia_binary} (modified as in \citep*[Fig.~3]{Baake_2019_renormalisation}).
    \item \href{https://tilings.math.uni-bielefeld.de/substitution/cap/}{CAP} \citep*{encyclopedia_cap}.
    \item \href{https://tilings.math.uni-bielefeld.de/substitution/watanabe-ito-soma-12-fold/}{Watanabe--Ito--Soma 12-fold} \citep*{encyclopedia_watanabe_ito_soma_12fold}.
    \item \href{https://tilings.math.uni-bielefeld.de/substitution/ammann-a3/}{Ammann A3} \citep*{encyclopedia_ammann_a3}.
\end{enumerate}
\end{introtheorem}
The Godr\`eche--Lan\c{c}on--Billard substitution rule is particularly interesting,
as it gives rise to hyperuniform point processes with singular continuous diffraction,
which makes it inaccessible to previous methods.

In the literature of aperiodic order,
it is customary to consider properties of tilings and point sets
which are invariant up to \emph{MLD equivalence}
\citep*[Def.~5.7]{Baake2013}.
However, in this article,
we prove that hyperuniformity is not an MLD invariant,
contradicting the standing conjecture among experts \citep*{invariance-conjecture}
(see \Cref{subsec:non-mld-invariance} for details).

\begin{introtheorem}\label{introtheorem:non-mld-equivalence}
There exist two disjoint FLC sets \(\Lambda_1\) and \(\Lambda_2\)
with the following properties:
\begin{itemize}
    \item \(\Lambda_1\) and \(\Lambda_1 \cup \Lambda_2\) are MLD equivalent.
    \item The point process associated to \(\Lambda_1\) is not hyperuniform.
    \item The point process associated to \(\Lambda_1 \cup \Lambda_2\) is hyperuniform.
\end{itemize}
\end{introtheorem}

The counterexample is constructed using a one-dimensional substitution rule,
similar to the ones considered by \citet*{Ouz2019},
whose substitution matrix has a small spectral gap.
\begin{figure}
	\centering
    \begin{subfigure}{0.32\linewidth}
        \adjincludegraphics[width=\linewidth,viewport=0 0 {0.4\height} {0.4\height},clip]{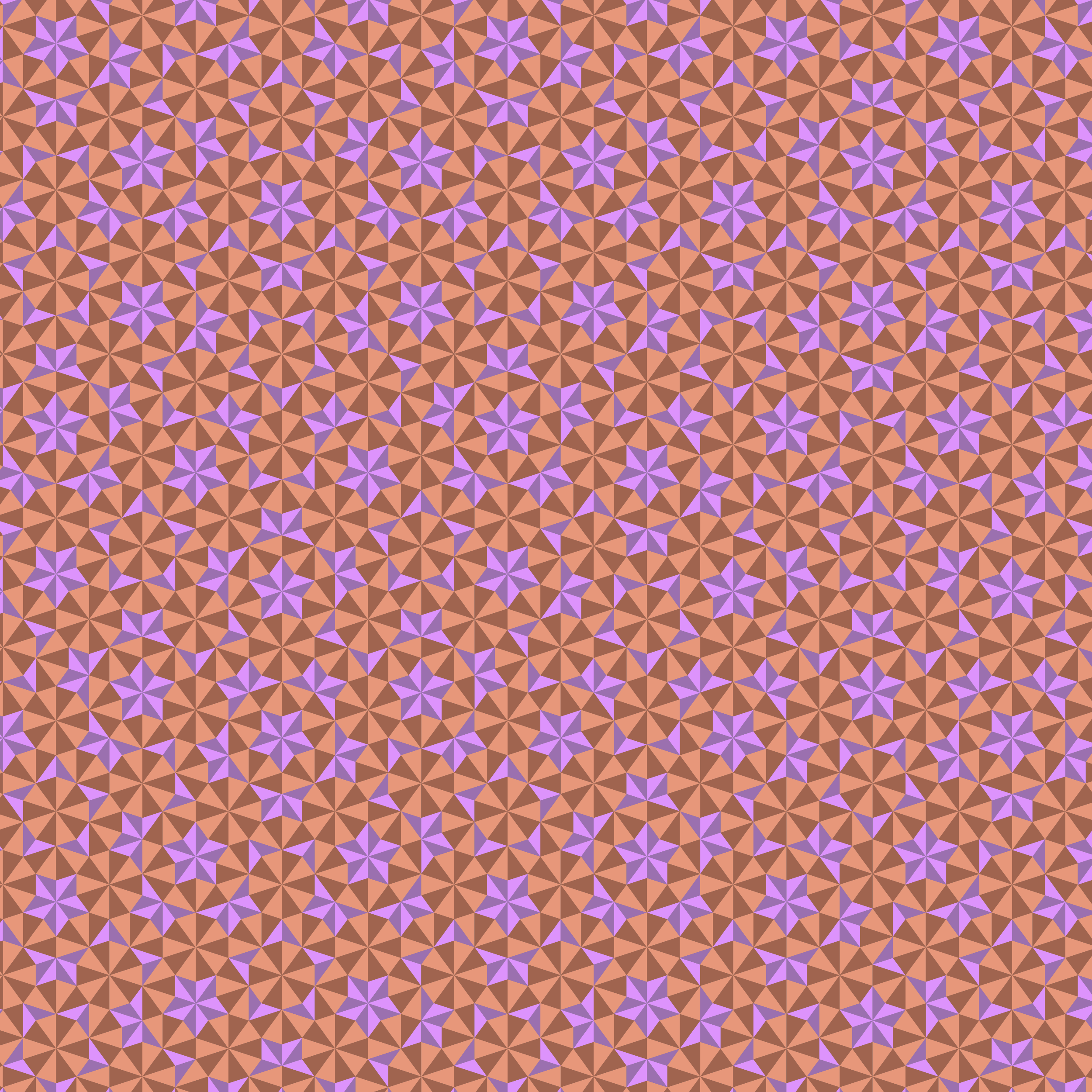}
        \caption{A Penrose tiling}
    \end{subfigure}
    \begin{subfigure}{0.32\linewidth}
        \adjincludegraphics[width=\linewidth,viewport=0 0 {0.15\height} {0.15\height},clip]{encyclopedia/danzer.pdf}
        \caption{A Danzer tiling}
    \end{subfigure}
    \begin{subfigure}{0.32\linewidth}
        \adjincludegraphics[width=\linewidth,viewport=0 0 {0.15\height} {0.15\height},clip]{encyclopedia/glb.pdf}
        \caption{A GLB tiling}
    \end{subfigure}
	\caption{The substitution tilings obtained from some of the substitution rules considered in \Cref{theorem:examples-intro}.
    Figures \href{https://tilings.math.uni-bielefeld.de/substitution/danzers-7-fold/}{(b)}
    and \href{https://tilings.math.uni-bielefeld.de/substitution/binary/}{(c)}
    are by \citet*{encyclopedia},
    licensed under \href{https://creativecommons.org/licenses/by-nc-sa/2.0/}{CC BY-NC-SA 2.0}.
    }\label{fig:example-thm-tilings-intro}
\end{figure}

\subsubsection{Self-similar point processes and substitution rules}
We focus on point processes coming from \emph{substitution rules}.
A substitution rule \(\cS = (\lambda, \bm{\Delta})\)
is characterized by a \emph{scaling constant} \(\lambda > 1\)
and a \emph{displacement matrix} \(\bm{\Delta} = (\Delta_{jk})_{j, k \in [\ell]}\),
where each \(\Delta_{jk}\) is a finite subset of \(\bR^d\).
We then define \emph{self-similar point processes},
which are a stochastic analogue of substitution Delone \(\ell\)-sets
studied by authors such as \citet*{lagarias2002substitutiondelonesets} and \citet*{LeeNonFLC}.

\begin{definition}
A \emph{self-similar point process with rule \(\cS\)}
is an \(\ell\)-tuple of stationary, locally square-integrable point processes
\(\bm{\Phi} = (\Phi_1, \ldots, \Phi_\ell)\)
such that \[
    \Phi_k \equidist \sum_{j=1}^\ell \sum_{x \in \Delta_{kj}} T_x D_\lambda \Phi_j, \quad \forall k \in [\ell]
\]
where \(\equidist\) denotes equality in distribution.
\end{definition}

If \(\cS\) is \emph{stone} and \emph{primitive} (to be defined later),
there exists a unique ergodic self-similar point process \(\bm{\Phi} = (\Phi_1, \ldots, \Phi_\ell)\) with rule \(\cS\),
which we call the \emph{vector point process associated to \(\cS\)}.
Furthermore, there exist compact sets \(\tau_1, \ldots, \tau_\ell \subset \bR^d\),
which we call the \emph{canonical prototiles} of \(\cS\),
such that the sets \(\Phi_1, \ldots, \Phi_\ell\) tile \(\bR^d\) with the prototiles \(\tau_1, \ldots, \tau_\ell\):
that is, the set \[
    \{T_x \tau_j \mid j \in [\ell], x \in \Phi_j\}
\]
is a random tiling of \(\bR^d\).

\begin{figure}[h!]
    \centering
    \begin{subfigure}{0.45\linewidth}
        \centering
        \adjincludegraphics[width=\linewidth]{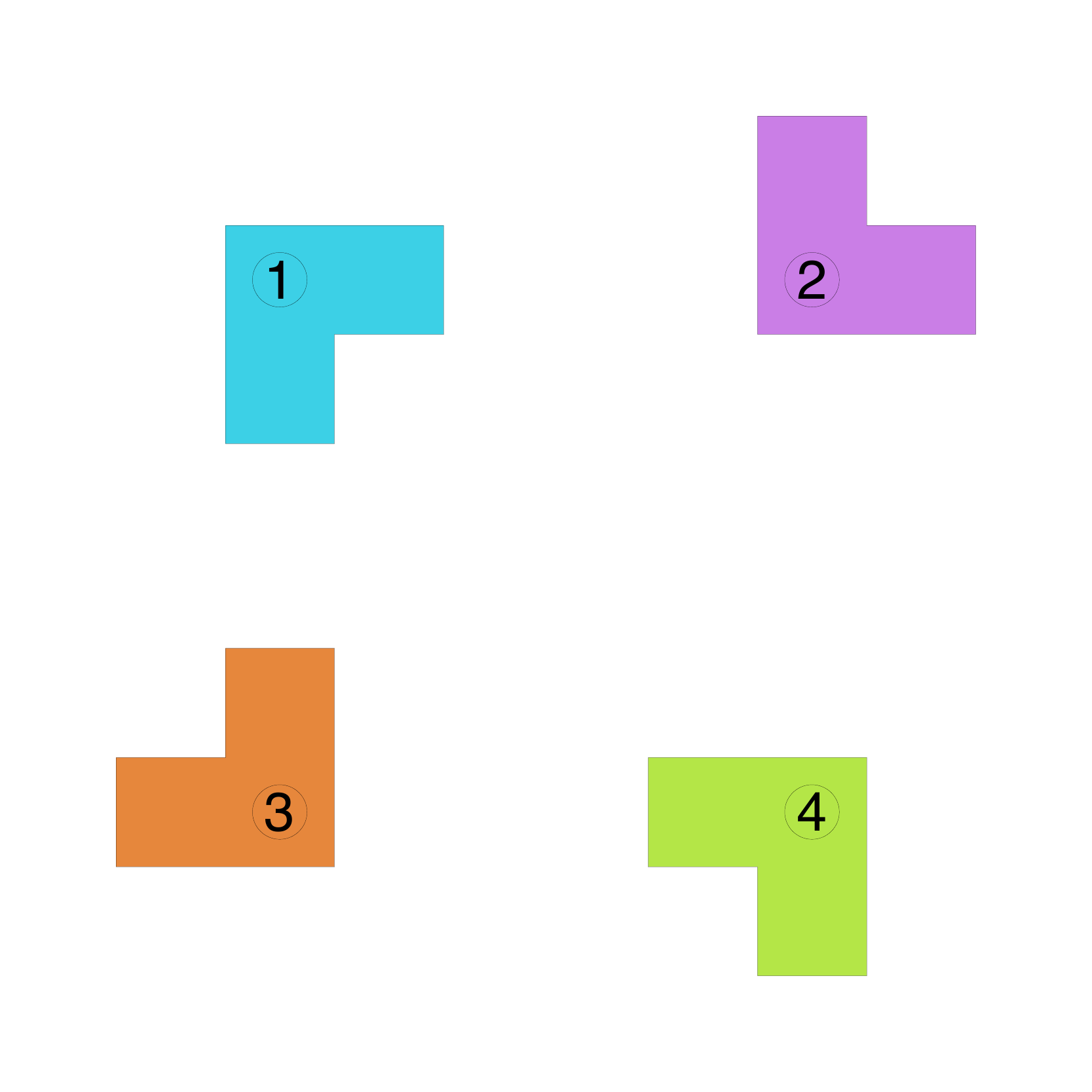}
        \caption{The four canonical prototiles of the chair substitution rule, with the origin marked.}
        \label{subfig:intro-chair}
    \end{subfigure}
    \begin{subfigure}{0.45\linewidth}
        \centering
        \adjincludegraphics[width=\linewidth, viewport=0 0 {0.4\height} {0.4\height},clip]{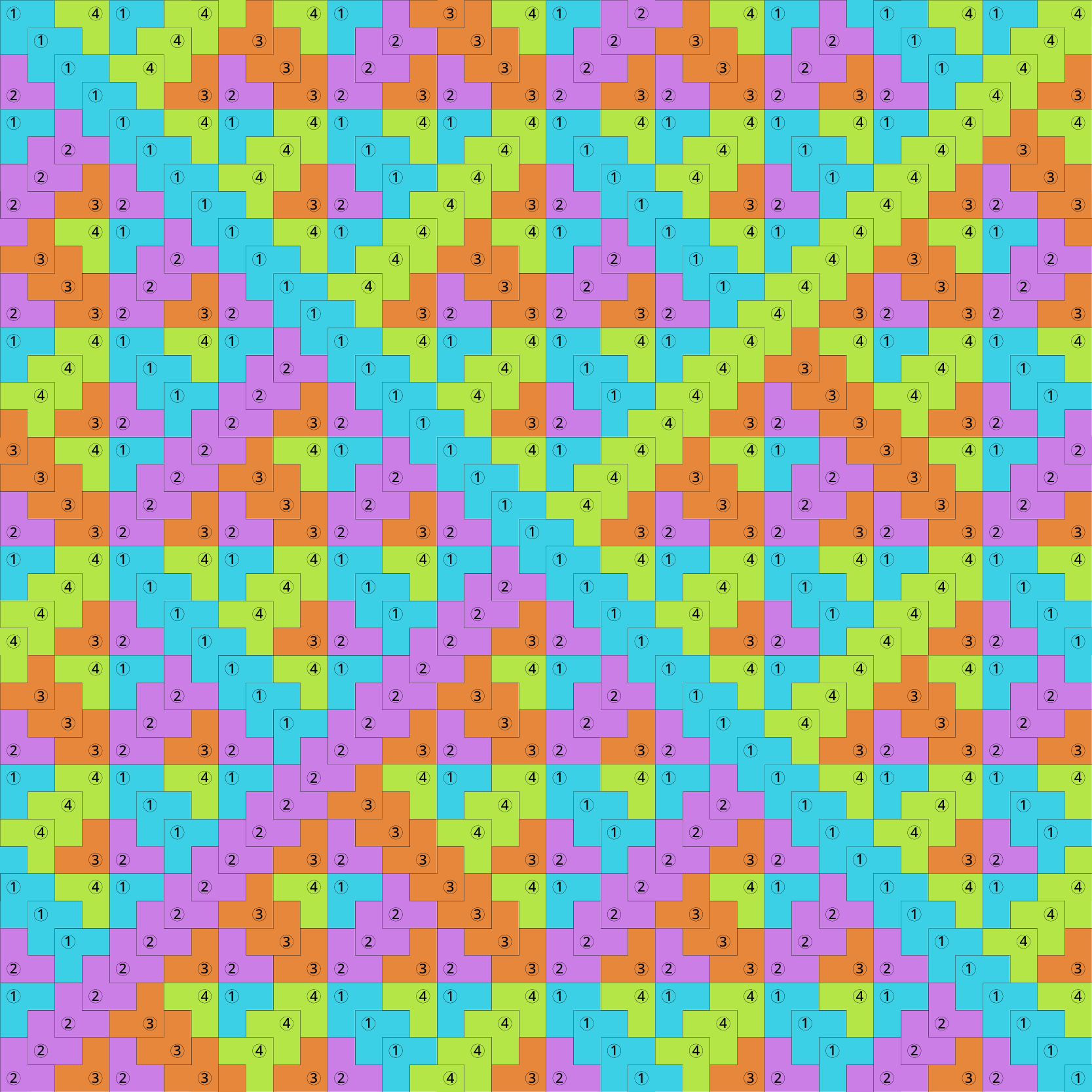}
        \caption{The four random sets \(\Phi_1, \Phi_2, \Phi_3, \Phi_4\) of the vector point process associated with the chair rule,
        which tile \(\bR^2\) with the four canonical prototiles as depicted.}
        \label{subfig:intro-chair-tiling}
    \end{subfigure}
    \caption{Canonical prototiles and vector point process for the chair substitution rule.}
    \label{fig:chair-rule-intro}
\end{figure}

\begin{example}
One of the most famous examples of a substitution rule is the \emph{chair rule},
a substitution rule on \(\bR^2\) with \(\ell = 4\),
and \(\lambda = 2\). The displacement matrix is given in \Cref{chair-rule-def}.
\Cref{fig:chair-rule-intro} depicts its canonical prototiles
and the associated vector point process.
\end{example}

\citet*{Baake_2019_origin}
proved a sufficient condition for the hyperuniformity of the random measure \({\Phi = \Phi_1 + \ldots + \Phi_\ell}\)
in terms of the \emph{substitution matrix} \(\Mfull \in \bR^{\ell \times \ell}\),
defined by \({(\Mfull)}_{jk} = \# \Delta_{jk}\).
By definition, \(\cS\) is \emph{primitive}
if \(\Mfull\) is primitive,
i.e. if there exists \(N \in \bN\) such that all entries of \(\Mfull^N\) are positive.
In this case, by the Perron--Frobenius theorem,
\(\Mfull\) has \(\lambda^d\) as its largest eigenvalue.

If \(\cS\) is primitive and stone,
the following holds.
\begin{theorem}[\citet*{Baake_2019_origin}]
Assume \(\cS\) has pure point diffraction.
Then the following holds:
if \(\abs{\mu} < \lambda^{\frac{d}2}\)
for all \(\mu \in \Spec \Mfull \setminus \{\lambda^d\}\),
then \(\Phi\) is hyperuniform.
\end{theorem}

Using this theorem,
Baake and Grimm were able to prove hyperuniformity
for several one-dimensional substitution rules,
such as the Fibonacci substitution rule.
However, when applied to two-dimensional examples
the spectral condition usually does not hold.
We prove a stronger criterion for hyperuniformity,
which is analogous to a condition for bounded displacement
due to \citet*{solomon-bd}.
For each eigenvalue \(\mu \in \Spec \Mfull\),
let \(E_\mu\) be the generalized eigenspace associated to \(\mu\).

\begin{introtheorem}
\label{thm:intro-hyper-criterion}
Assume \(d \in \{1,2\}\).
Then the following holds:
if \(\abs{\mu} < \lambda^{\frac{d}{2}}\)
for all \(\mu \in {\Spec \Mfull \setminus \{\lambda^d\}}\)
such that \((1, \ldots, 1)^\top \not\in E_\mu^\perp\),
then \(\Phi\) is hyperuniform.
\end{introtheorem}

This improves on the previous criterion
in two crucial ways:
\begin{itemize}
    \item Instead of requiring \(\abs{\mu} < \lambda^{\frac{d}{2}}\)
    for all \(\mu \in {\Spec \Mfull \setminus \{\lambda^d\}}\),
    it suffices to check this condition only for a subset of those eigenvalues.
    \item We do not need to assume pure point diffraction: this allows us to apply the criterion
    to new examples such as the Godr\`eche--Lan\c{c}on--Billard substitution rule,
    and makes it significantly easier to check the conditions.
\end{itemize}
Almost all known examples of substitution rules are in either one or two dimensions,
so the assumption that \(d \in \{1,2\}\) is not a significant restriction in practice.
It may be possible to adapt our methods to the \(d > 2\) case by making stronger assumptions on the spectrum of \(\Mfull\).

Most substitution rules in two dimensions have rotational and/or reflection symmetry.
In this case,
we can consider their \emph{spherical substitution matrix} \(\Msph\),
which is obtained by counting prototiles only up to the action of the group of symmetries.
This matrix is usually significantly smaller than \(\Mfull\).
Then \Cref{thm:intro-hyper-criterion}
implies the following simpler criterion,
which is powerful enough to apply to all examples in \Cref{theorem:examples-intro}.
\begin{introtheorem}
\label{thm:intro-hyper-criterion-symm}
Assume \(d \in \{1,2\}\).
Then
the following holds:
if \(\abs{\mu} < \lambda^{\frac{d}{2}}\)
for all \(\mu \in {\Spec \Msph \setminus \{\lambda^d\}}\),
then \(\Phi\) is hyperuniform.
\end{introtheorem}

\subsubsection{Diffraction}
If \(\Phi\) is a locally square-integrable stationary random measure on \(\bR^d\),
we can define its \emph{diffraction} or \emph{Bartlett spectral measure} \(\hat{\eta}\)
as the unique Radon measure on \(\bR^d\)
satisfying \[
    \hat{\eta}(\abs{\hat{f}}^2) = \Var (\Phi(f)), \quad \forall f \in C_c(\bR^d).
\]
In our setting, the diffraction can also be defined as the Fourier transform of the \emph{autocorrelation}
of any sample, which is defined using an ergodic average due to \citet*{Hof1995}
(note that, in the literature of aperiodic order, a different convention for the diffraction is used,
which differs from ours only by a Dirac mass at the origin).

\citet*{hartnick2022hyperuniformity} proved the following criterion for hyperuniformity
in terms of the diffraction measure.
\begin{theorem}[{\citep*{hartnick2022hyperuniformity}}]
Let \(\Phi\) be a locally square-integrable stationary point process on \(\bR^d\),
and let \(\hat{\eta}\) be its diffraction measure.
Then \(\Phi\) is hyperuniform if and only if \(\lim_{r \to 0}\frac{\hat{\eta}(B_r)}{r^d} = 0\).
\end{theorem}

Let \(\cS\) be a substitution rule in \(\bR^d\)
with scaling constant \(\lambda\) and substitution matrix \(\Mfull\).
Define the constants:
\begin{align*}
    \beta_\parallel &:= d - \max\left\{ \log_\lambda \abs{\mu} \;\middle|\; \mu \in \Spec \Mfull \setminus \{\lambda^d\},\; (1, \ldots, 1)^\top \not\in E_\mu^\perp \right\}, \\
    \beta_\perp &:= d+1 - \max\left\{ \log_\lambda \abs{\mu} \;\middle|\; \mu \in \Spec \Mfull \setminus \{\lambda^d\},\; (1, \ldots, 1)^\top \in E_\mu^\perp \right\}, \\
    \beta &:= \min\left(\beta_\perp, \beta_\parallel\right).
\end{align*}

We prove the following bound on \(\hat{\eta}\) around the origin
(see \Cref{sec:lyapunov-asymptotic-notation} for the definition of the symbol \(\lessapprox\)).
\begin{introtheorem}\label{intro-origin-bound}
Let \(\bm{\Phi} = (\Phi_1, \ldots, \Phi_\ell)\) be a self-similar point process with rule \(\cS\),
and \(\hat{\eta}\) be the diffraction of \(\Phi = \Phi_1 + \ldots + \Phi_\ell\).
Then we have \[
    \hat{\eta}(B_r) \lessapprox r^{2 \beta}
\]
as \(r \to 0\).
\end{introtheorem}
\Cref{thm:intro-hyper-criterion},
as well as \Cref{theorem:examples-intro},
follow as corollaries of \Cref{intro-origin-bound}.
Furthermore, this bound also yields Class I hyperuniformity in many cases.
Beyond determining when a process is or is not hyperuniform,
the behaviour of the diffraction around the origin encodes interesting information about the point process,
such as {rigidity phenomena}
\citep*{Ghosh_2017, lachiezerey2025rigidityrandomstationarymeasures}.
Physicists have studied this topic extensively,
heuristically and with numerical experiments
\citep*{Torquato_2018, Ouz2019, torquato-penrose-numeric, Fuchs2019}.

\subsubsection{Renormalisation}
Let \(\cS\) be a primitive, stone substitution rule
on \(\bR^d\)
with scaling constant \(\lambda\),
displacement matrix \(\bm{\Delta} = (\Delta_{jk})_{j, k \in [\ell]}\)
and full substitution matrix \(\Mfull\).
Let \(\bm{\Phi} = (\Phi_1, \ldots, \Phi_\ell)\)
be a self-similar point process with rule \(\cS\).
Then the \emph{diffraction matrix} \(\hat{\bm{H}} = (\hat{H}_{jk})_{j,k \in [\ell]}\)
of \(\bm{\Phi}\)
is defined by the formula \[
    \hat{H}_{jk}(\abs{\hat{f}}^2) = \Cov(\Phi_j(f), \Phi_k(f)), \quad \text{for all \(f \in C_c(\bR^d)\), \(j, k \in [\ell]\)}.
\]

Due to the self-similarity of \(\bm{\Phi}\),
the diffraction matrix satisfies the following \emph{renormalisation relations}.

\begin{theorem}[\citep*{Baake_2019_renormalisation}]
\label{thm:diffraction-renormalization-intro}
Let \(\cS\) be a primitive, stone substitution rule in \(\bR^d\)
with \(\ell\) prototiles,
\(\lambda\) its scaling constant,
and \(\Mfull\) its substitution matrix.
Let \(\hat{\bm{H}}\) be its diffraction matrix.
Then there exists a smooth, matrix-valued function \(\bm{A} : \bR^d \to \bC^{\ell \times \ell}\)
which satisfies \({\bm{A}(0)=\frac{1}{\lambda^d}\Mfull}\),
called the \emph{normalized Fourier matrix} of the substitution rule,
such that the following holds:
\[
    \hat{\bm{H}} \;=\; \bm{A}\bigl(D_{\lambda^{-1}}\hat{\bm{H}}\bigr)\,\bm{A}^*.
\]
\end{theorem}

\citet*{Baake_2019_renormalisation} used the above renormalisation relation
to study the pure point part and the absolutely continuous part of the diffraction.
In these cases, there is a natural way to write \(\hat{\bm{H}}\)
using a density function
(with respect to either the counting measure or the Lebesgue measure, respectively),
which then satisfies a recurrence relation.

One of the core innovations of this article
is that we extend this method to the singular continuous part of the diffraction.
We do this by introducing a \emph{self-similar density for \(\hat{\bm{H}}\)}:
this is a density function for the diffraction matrix
which is compatible with the renormalisation relation
and exists without making assumptions on the spectral type of the diffraction.

\subsubsection{Linear cocycles}
To analyze diffraction near the origin,
we study iterated products of the normalized Fourier matrix \(\bm{A}\),
namely the \emph{Fourier cocycle} \[\bm{A}^{(N)}(\xi) = \bm{A}(\lambda^{-N+1}\xi)\bm{A}(\lambda^{-N+2}\xi)\cdots \bm{A}(\xi)\]
for \(\xi \in \bR^d\) and \(N \in \bN\).
This is similar to the \emph{internal cocycle} considered
by \citet*{fibonacci-diffraction}
and the \emph{spectral cocycle}
considered by \citet*{solomyak2022spectral}.
To prove \Cref{thm:intro-hyper-criterion},
we study the asymptotics of \(\langle \bm{A}^{(N)}(\xi)x, w\rangle\)
as \(N \to \infty\).

The normalized Fourier matrix \(\bm{A}\)
is smooth; therefore \(\bm{A}(\xi) = \bm{A}(0) + O(\norm{\xi})\) as \(\xi \to 0\).
Then, if \(N\) is large,
one would intuitively expect
\(\bm{A}^{(N)}(\xi)\) to be close to \(\bm{A}(0)^N\),
as most of the matrices in the product
are close to \(\bm{A}(0)\).

To formalize this heuristic
we prove an Oseledets-like theorem for products of exponentially converging matrices,
similar to a result due to \citet*{Dubiner1991}.
Applied to the Fourier cocycle,
it yields \Cref{intro-origin-bound}.
This result holds beyond the specific setting of linear cocycles arising from substitution rules,
hence it may be of independent interest.
Below we give a simplified version:
see \Cref{thm:main-theorem-lyapunov} for the full statement,
and \Cref{sec:lyapunov-asymptotic-notation}
for the definition of the asymptotic symbols \(\lessapprox\) and \(\approx\).

\begin{introtheorem}
    \label{lyapunov-intro}
Let \(\chi_1 > \cdots > \chi_l\)
be the distinct values of \(\{\log \abs{\mu} \mid \mu \in \Spec \bm{A}(0), \mu \neq 0\}\)
and let \(\chi_{l+1} = -\infty\).
Let \(E_j := \oplus \{E_\mu \mid \mu \in \Spec \bm{A}(0), \log \abs{\mu} = \chi_j\}\).
There exist uniquely defined idempotent operators \(P_j : \bC^\ell \to \bC^\ell\)
for \(j \in [l+1]\)
such that \(\im P_j = E_j\) and \(\sum_{j=1}^{l+1} P_j = I\).

Then there exists a radius \(R > 0\)
such that, for all \(\xi \in B_R\),
there exist linear maps \({P_j(\xi) : \bC^\ell \to \bC^\ell}\) for all \(j \in [l+1]\)
such that \(\sum_{j=1}^{l+1} P_j(\xi) = I\)
and the following asymptotic inequalities hold:
\begin{enumerate}
    \item For all \(j \in [l+1]\) we have
    \begin{equation*}
             \norm{\bm{A}^{(N)}(\xi)P_j(\xi)x} \approx e^{\chi_j N}\norm{P_j(\xi)x}
    \end{equation*}
    as \(N \to \infty\).
    \item For all \(j, k \in [l+1]\) we have
    \begin{align*}
            {\norm{P_k\bm{A}^{(N)}(\xi)P_j(\xi)x}} &\lesssim \lambda^{-N}\norm{\bm{A}^{(N)}(\xi)P_j(\xi)x}
    & \mbox{if \(j > k\)}\\
    {\norm{P_k\bm{A}^{(N)}(\xi)P_j(\xi)x}} &\lessapprox \max(e^{\chi_k - \chi_j}, \lambda^{-1})^N\norm{\bm{A}^{(N)}(\xi)P_j(\xi)x}
        & \mbox{if \(j < k\)}
    \end{align*}
    as \(N \to \infty\).
\end{enumerate}
\end{introtheorem}

\subsubsection{Acknowledgements and related work}
This article summarizes and improves upon the author's doctoral thesis \citep*{thesis}.
I would like to thank my PhD advisor, Tobias Hartnick, for sparking my interest in aperiodic order and helpful mathematical discussions.
This work was partially supported by the German Research Council (Deutsche Forschungsgemeinschaft, DFG) under RTG 2229 (``Asymptotic Invariants and Limits of Groups and Spaces'').

\section{Preliminaries}\label{chap:preliminaries}

\subsection{Notation and basics}
\label{sec:notation}

In this article, we write \(\bN\) for the natural numbers without \(0\) and \(\bN_0\) for the natural numbers including \(0\).
We write \(\bR\) for the real numbers and \(\bR_{\geq 0}\) for the non-negative reals.
For \(n \in \bN\) we let \([n] = \{1, \ldots, n\}\).
We denote the number of elements in a finite set \(S\) by \(\# S\).
If \(X\) and \(Y\) are two random variables with values in the measurable space \(S\),
we write \(X \equidist Y\) if they are equidistributed,
i.e. they have the same probability distribution measure.
Our inner products are always conjugate-linear in the second argument.
We denote the conjugate-transpose of a matrix \(A \in \bC^{\ell_1 \times \ell_2}\) by \(A^*\).

We let \(B_R \subset \bR^d\) be the closed ball of radius \(R\) around the origin,
and \(B_R^\times = B_R \setminus \{0\}\) be the punctured ball.
We let \(m_{\bR^d}\) be the Lebesgue measure on \(\bR^d\).

For \(x \in \bR^d\) and \(\lambda > 0\), we write \(T_x\) for the translation operator on \(\bR^d\)
and \(D_\lambda\) for the dilation operator: these can also be applied to subsets of \(\bR^d\) as usual.
If \(f\) is a function on \(\bR^d\),
we define functions \(T_x f\) and \(D_\lambda f\)
by \((T_x f)(y) = f(y - x)\) and \((D_\lambda f)(y) = f(\lambda^{-1} y)\) respectively.
If \(f\) is a complex-valued function on \(\bR^d\),
we denote its complex conjugate by \(\overline{f}\).
We also define the function \(\tilde{f}\)
by \(\tilde{f}(x) = \overline{f(-x)}\).
If \(f \in L^1(\bR^d)\)
we denote its Fourier transform by \(\hat{f}\) or \(\cF(f)\),
defined by \(\hat{f}(\xi) = \int_{\bR^d} f(x) e^{-2\pi i \langle x, \xi \rangle} \,\mathrm{d}x\).
We denote its inverse Fourier transform by \(\check{f}\) or \(\cF^{-1}(f)\).
We also define the \emph{convolution} of two functions \(f, g \in L^1(\bR^d)\)
by \(
    (f * g)(x) = \int_{\bR^d} f(y) g(x - y) \,\mathrm{d}y
    \).
\section{Diffraction of random vector measures}\label{chap:diffraction}
In this section,
we introduce the basics of stationary random vector measures
and their diffraction.
This is a straightforward generalization of the classical theory of random measures and point processes
\citep*{poisson,Brmaud2020}.
For a more detailed treatment,
including proofs of the statements presented here,
see \citep*[Chapter 3]{thesis}.

\subsection{Complex and vector measures}
\label{sec:vector-measures}

\index[defs]{measure!complex}
Let \(C_c(X)\) be the space of (complex-valued) continuous functions on \(X\)
with compact support.
A \emph{complex measure \(\mu\) on \(\bR^d\)}
(in the sense of \citet*{bourbaki-integration}) is, by definition,
a continuous linear functional \(\mu : C_c(\bR^d) \to \bC\):
any such map can be uniquely extended to bounded measurable functions with bounded support,
such that \(\mu(A) := \mu(\bone_A)\) is well-defined for all bounded Borel sets \(A \subset \bR^d\).
We let \(\cM(X)\) be the space of complex measures on \(X\)
and \(\cM_+(X)\) be the space of positive measures on \(X\).
We equip \(\cM(X)\) with the weak-* topology.

\begin{remark}
This is the way \citet*{bourbaki-integration} defines a measure,
and it is ubiquitous in the fields of harmonic analysis and aperiodic order \citep*{Hof1995,Berg1975,moody_strungaru_2017}.
Note that, using this definition,
the map \(\mu : C_c(\bR) \to \bC\)
given by \(\phi \mapsto \int_{-\infty}^\infty \phi(x) \sin(x)\,\mathrm{d}x\)
is a complex measure,
but \(\mu(\bR)\) is not well defined:
\(\mu(A)\) is only necessarily defined when \(A\) is a \emph{bounded} Borel set.
\end{remark}

\index[defs]{measure!vector}
\index[defs]{measure!matrix}
A \emph{\(\bC^\ell\)-measure \(\bm{\mu}\) on \(X\)}, or \emph{vector measure},
is a vector of complex measures \(\bm{\mu} = (\mu_j)_{j \in [\ell]}\).
Equivalently this is a continuous linear map \(\bm{\mu} : C_c(X) \to \bC^\ell\).
We denote by \(\cM(X, \bC^\ell) := \cM(X)^\ell\) the space of \(\bC^\ell\)-measures on \(X\),
equipped with the weak-* topology on the components.
Matrix measures are defined analogously.

\index[defs]{total variation measure}
Given a vector measure \(\bm{\mu}\),
we can define its \emph{total variation measure} \(\norm{\bm{\mu}}\)
which is the least positive measure
satisfying \(\norm{\bm{\mu}(f)} \leq \norm{\bm{\mu}}(\abs{f})\).
See \citet*{bourbaki-integration} for details.
Then we can uniquely extend \(\bm{\mu}\) to \(L^1(\norm{\bm{\mu}})\) by continuity.
Appropriate extensions of the Fubini and Radon--Nikodym theorems
hold in this setting \citep*{bourbaki-integration}.

\begin{definition}
Let \(\bm{\mu} = (\mu_j)_{j=1}^\ell\) be a \(\bC^\ell\)-measure on \(\bR^d\)
and \(w \in \bC^\ell\) be a weight vector.
We define the \emph{weighted measure} \(\langle \bm{\mu}, w\rangle\)
by \[
    \langle \bm{\mu}, w\rangle (f) = \sum_{j=1}^\ell \overline{w}_j \mu_j(f), \quad \forall f \in C_c(\bR^d).
\]
\end{definition}
We also use vector/matrix notation to write products between matrix measures and vector measures.

Let \(\mu\) be a complex measure on \(X\).
In analogy with the function case,
for \(x \in \bR^d\) and \(\lambda > 0\),
we define \(T_x \mu \) and \(D_\lambda \mu\)
by \((T_x \mu)(f) = \mu(T_{-x} f)\)
and \((D_\lambda \mu)(f) = \mu(D_{\lambda^{-1}} f)\)
for all \(f \in C_c(\bR^d)\) respectively.
We also define the conjugate and adjoint
by \(\overline{\mu}(f) = \overline{\mu(\overline{f})}\)
and \(\tilde{\mu}(f) = \overline{\mu(\tilde{f})}\)
for all \(f \in C_c(\bR^d)\) respectively.
In addition, if \(\mu\) is a translation bounded complex measure on \(\bR^d\),
we denote its Fourier transform by \(\hat{\mu}\),
assuming it exists (see \citet*{Berg1975}).
These operations are defined componentwise for vector and matrix measures as well.

\subsection{Stationary random vector measures}
\label{sec:stationary-random-vector-measures}
For the rest of this section let \(d, \ell \in \bN\).

A random \(\bC^\ell\)-measure \(\bm{\Phi}\) is \begin{itemize}
    \item \index[defs]{stationary!random measure}
 \emph{stationary}
if,
for all \(x \in \bR^d\),
\(T_x \bm{\Phi} \equidist \bm{\Phi}\); and
\item \index[defs]{locally square-integrable} \emph{locally square-integrable}
if, for all \(f \in C_c(\bR^d)\), we have \(\bE[\norm{\bm{\Phi}(f)}^2] < \infty\).
\end{itemize}
\textbf{From now on we assume all random (vector) measures are stationary and locally square-integrable.}
Furthermore, we say \(\bm{\Phi}\) is an \emph{\(\ell\)-point process} if every component is a simple, locally finite point measure:
in this case, one can alternatively interpret \(\bm{\Phi}\)
as a random vector of locally finite subsets of \(\bR^d\).

\begin{example}
Let \(\Gamma = \bZ^2\) and \(\Omega = \{T_x \Gamma \mid x \in \bR^2\}\).
We have \(\Omega = \bR^2/\bZ^2 = \bT^2\):
using the uniform distribution on \(\Omega\),
we can define a stationary random measure by \[
    \Phi : \Omega \to \cM(\bR^2), \quad \Gamma' \mapsto {\Phi}_{\Gamma'} = \sum_{y \in \Gamma'} \delta_y
\]
\end{example}

\begin{example}
The \emph{homogeneous Poisson process on \(\bR^d\) with intensity \(\lambda\)} \citep*{poisson}
is the unique stationary point process with the following two properties:
\begin{enumerate}
    \item The expectation of \({\Phi}\) is given by \(\bE[{\Phi}(f)] = \lambda m_{\bR^d}(f)\)
    \item For any \(f, f' \in C_c(\bR^d)\) such that \(f\) and \(f'\) have disjoint support,
    \({\Phi}(f)\) and \({\Phi}(f')\) are independent random variables.
\end{enumerate}
\end{example}

One of the most important objects associated to a stationary random measure
is its \emph{diffraction}, also called \emph{spectral Bartlett measure}
or \emph{structure factor} \citep*{Brmaud2020,coste_2021,bjorklund2024hyperuniformityrandommeasureseuclidean}.

\begin{definition}
Let \(\Phi\) be a random measure on \(\bR^d\).
Its \emph{(centered) diffraction measure} \(\hat{\eta}\)
is the unique measure which satisfies \[
    \hat{\eta}(\abs{\hat{f}}^2) = \Var(\Phi(f)), \quad \forall f \in C_c(\bR^d).
\]
\end{definition}

\begin{remark}\label{rem:hof-diffraction}
The above quantity is related to the diffraction considered in the literature of aperiodic order as follows.
Let \(\omega\) be a complex measure
and, for \(R > 0\),
write \(\omega_R := \omega|_{B_R}\).
Then its \emph{(uncentered) autocorrelation measure \(\gamma\)}
is defined by the limit \[
    \gamma := \lim_{R \to \infty} \frac{\omega_R * \omega_R^*}{\Vol(B_R)}
\]
in the vague topology,
assuming this limit exists and is unique.
If \(\gamma\) is well-defined,
it is a positive definite measure,
hence we can define its \emph{Hof diffraction measure \(\hat{\gamma}\)}
as the Fourier transform of \(\gamma\).
This definition is due to \citet*{Hof1995}
and is the classical definition of diffraction in the literature of aperiodic order.
As pointed out in \citep*{BaakeBirknerGrimm2015,BaakeBirknerMoody2009},
if \(\Phi\) is an \emph{ergodic} stationary random measure,
then its autocorrelation \(\gamma\) is almost surely well-defined
and satisfies \[
    \hat{\gamma}(\abs{\hat{f}}^2) = \bE[\abs{\Phi(f)}^2] , \quad \forall f \in C_c(\bR^d).
\]
The Hof diffraction and the spectral Bartlett measure
are essentially the same, except the Hof diffraction is defined using the squared expectation
and the spectral Bartlett measure is defined using the variance.
If we let \(\iota\) be the \emph{intensity} of \(\Phi\), i.e. \(\iota = \bE[\Phi(B_1)]\),
they are related by the formula \(\hat{\eta} = \hat{\gamma} - \abs{\iota}^2 \delta_0\),
so they coincide up to the removal of an atom at the origin.
\end{remark}

For a random vector measure,
we can define the following generalization of the diffraction measure,
which takes the different components into account.

\begin{definition}\label{def:diffraction}
Let \(\bm{\Phi}\) be a stationary random \(\bC^\ell\)-measure on \(\bR^d\).
Its \emph{diffraction matrix}
is the unique matrix measure \(\hat{\bm{H}} = (\hat{H}_{jk})_{j,k=1}^\ell\)
which satisfies \[
    \hat{H}_{jk}(\abs{\hat{f}}^2) = \Cov({\Phi}_j(f), {\Phi}_k(f)), \quad \forall f \in C_c(\bR^d), \; j,k \in [\ell].
\]
\end{definition}

If \(\bm{\Phi}= (\Phi_j)_{j=1}^\ell\) is a stationary random \(\bC^\ell\)-measure,
for every \(w \in \bC^\ell\) we obtain a stationary random measure
\(\langle \bm{\Phi}, w\rangle = \sum_{j=1}^\ell \overline{w}_j \Phi_j\).
For example, if \(w = (1, \ldots, 1)^\top\) is the constant vector,
then \(\langle \bm{\Phi},w\rangle = \sum_{j=1}^\ell \Phi_j\),
while if \(w = e_j\) for some \(j \in [\ell]\),
then \(\langle \bm{\Phi},w\rangle = \Phi_j\).
The matrix diffraction \(\hat{\bm{H}}\) is related to the diffraction measures as follows:
for all \(w \in \bC^\ell\),
\(\langle \hat{\bm{H}}w,w\rangle\) is the diffraction measure of \(\langle \bm{\Phi}, w\rangle\).
In particular, \(\langle \hat{\bm{H}}w,w\rangle\) is a positive measure for all \(w \in \bC^\ell\).

\subsection{Hyperuniformity}

Now we define the notion of \emph{hyperuniformity},
introduced by \citet*{Torquato_2003} in the context of point processes.
Intuitively, hyperuniformity indicates a certain degree of \enquote{order},
in the sense that the variance of the measure on large sets
is less than one would expect from random chance.

We define hyperuniformity for complex-valued measures first,
which is the classical definition,
and then we extend it to vector-valued measures by using weights.
\begin{definition}[\citep*{Torquato_2003}]
\index[defs]{hyperuniform}
Let \(\Phi\) be a stationary random (complex) measure on \(\bR^d\).
\begin{enumerate}
    \item We say \(\Phi\) is \emph{hyperuniform} if \[
    \lim_{R \to \infty} \frac{\Var \Phi(B_R)}{R^d} = 0
    \]
    \item We say \(\Phi\) is \emph{Class I hyperuniform} if \[
    \Var \Phi(B_R) = O(R^{d-1})
    \] as \(R \to \infty\).
\end{enumerate}

If \(\bm{\Phi}\) is a stationary random \(\bC^\ell\)-measure on \(\bR^d\)
and \(w \in \bC^\ell\),
we say \(\bm{\Phi}\) is \emph{hyperuniform for weights \(w\)}
if \(\langle \bm{\Phi},w\rangle\) is hyperuniform:
if \(w = (1, \ldots, 1)^\top\) we say \(\bm{\Phi}\) is \emph{hyperuniform for constant weights}.
We define Class I hyperuniformity for vector measures analogously.
\end{definition}

We prove hyperuniformity using the following criterion
due to \citet*{hartnick2022hyperuniformity},
which relates hyperuniformity to the behaviour of the diffraction measure at the origin.
\begin{theorem}[\citep*{hartnick2022hyperuniformity}]
Let \(\Phi\) be a stationary random measure on \(\bR^d\) and
\(\hat{\eta}\) its diffraction measure.
\begin{enumerate}
    \item For \(\alpha \in [0,1]\), we have
    \(\Var(\Phi(B_R)) = O(R^{d-\alpha})\) as \(R \to \infty\)
    if and only if
    \(\hat{\eta}(B_r) = O(r^{d+\alpha})\) as \(r \to 0\).
    \item For \(\alpha \in [0, 1)\), we have
    \(\Var(\Phi(B_R)) = o(R^{d-\alpha})\) as \(R \to \infty\)
    if and only if
    \(\hat{\eta}(B_r) = o(r^{d+\alpha})\) as \(r \to 0\).
    \item \(\Phi\) is hyperuniform
    if and only if \(\hat{\eta}(B_r) = o(r^d)\) as \(r \to 0\).
    \item \(\Phi\) is Class I hyperuniform
    if and only if \(\hat{\eta}(B_r) = O(r^{d+1})\) as \(r \to 0\).
\end{enumerate}
\end{theorem}

\section{Self-similar point processes and renormalisation relations}\label{chap:substitutions}
In this section we introduce \emph{self-similar point processes},
which are point processes arising from the theory of substitution tilings.
We see that their diffraction matrices satisfy \emph{renormalisation relations},
which form the basis for our approach to their hyperuniformity.
For all of this section, let \(d, \ell \in \bN\).

\subsection{Self-similar point processes}
Recall that a \emph{substitution rule} \(\cS\)
is a pair \((\lambda, \bm{\Delta})\)
where \(\lambda > 1\) is the \emph{scaling constant}
and \(\bm{\Delta} \in \cP_{\mathrm{fin}}(\bR^d)^{\ell \times \ell}\)
is the \emph{displacement matrix}.
We also define the \emph{(full) substitution matrix} \(\Mfull \in \bN_0^{\ell \times \ell}\)
by \(M_{jk} = \# \Delta_{jk}\).

Then a \emph{self-similar point process with rule \(\cS\)}
is a stationary random \(\bC^\ell\)-measure
\(\bm{\Phi} = (\Phi_1, \ldots, \Phi_\ell)\) 
such that \[
    \Phi_k \equidist \sum_{j=1}^\ell \sum_{x \in \Delta_{kj}} T_x D_\lambda \Phi_j, \quad \forall k \in [\ell].
\]
In \Cref{sec:construction},
we will see that self-similar point processes always exist,
assuming some mild conditions on \(\cS\).

If \(\bm{\Phi}\) is a self-similar point process,
its diffraction matrix satisfies \emph{renormalisation relations}
as originally introduced by \citet*{Baake_2019_renormalisation}.

\begin{definition}
\index[defs]{normalized Fourier matrix}
The \emph{normalized Fourier matrix} of \(\cS\)
is the matrix function \(\bm{A} : \bR^d \to \bC^{\ell \times \ell}\)
with entries \(A_{jk}(\xi)\)
given by \[
    A_{jk}(\xi) = \lambda^{-d}\sum_{x \in \Delta_{jk}} e^{2\pi i \langle x, \xi\rangle}.
\]
\end{definition}
By definition, the normalized Fourier matrix \(\bm{A}(\xi)\)
is a componentwise analytic function of \(\xi\)
and satisfies \(\bm{A}(0) = \lambda^{-d} \Mfull\).

\begin{theorem}\label{thm:diffraction-renormalisation-measure}
The diffraction matrix measure \(\hat{\bm{H}}\) satisfies
\[
   \hat{\bm{H}} = \bm{A} (D_{\lambda^{-1}}\hat{\bm{H}}) \bm{A}^*.
\]
\end{theorem}

\begin{proof}[Proof of \Cref{thm:diffraction-renormalisation-measure}]
It suffices to check the equation holds for all functions
of the form \(\abs{\hat{f}}^2\) for \(f \in C_c(\bR^d)\),
as this uniquely determines the Fourier transform of a measure.
Coefficient-wise, the equation we want to prove is
\begin{equation}\label{eq:renormalisation-diffraction-coefficients}
    \hat{H}_{mn}(\abs{\hat{f}}^2) = \sum_{j,k=1}^\ell A_{mj} A_{nk}^* (D_{\lambda^{-1}}\hat{H}_{jk})(\abs{\hat{f}}^2), \quad \forall m,n \in [\ell].
\end{equation}

So let \(f \in C_c(\bR^d)\).
Recall that we have \[
    \hat{H}_{mn}(\abs{\hat{f}}^2) = \Cov(\Phi_m(f),\Phi_n(f))
\]
Then the result follows by applying the definition of a self-similar point process
to both \(\Phi_m(f)\) and \(\Phi_n(f)\),
simplifying using the properties of the Fourier transform for functions,
and writing the result in terms of \(\hat{\bm{H}}\).
\end{proof}
\subsection{Renormalisation via density functions}
\label{sec:renormalisation-density}
Now
we define \emph{self-similar densities}
as a tool to study the diffraction measure \(\hat{\bm{H}}\).

\begin{definition}
\index[defs]{density}
Let \(\bm{M} \in \cM(\bR^d, \bC^{\ell \times \ell})\) be a matrix valued measure on \(\bR^d\).
A \emph{density} of \(\bm{M}\) is a pair \((\bm{h}, \nu)\)
where
\begin{enumerate}
    \item \(\nu\) is a positive (scalar, and not necessarily \(\sigma\)-finite)
    measure on \(\bR^d\), and
    \item \(\bm{h}\)
    is a locally integrable \(\bC^{\ell \times \ell}\)-valued function on \(\bR^d\)
    such that \(\bm{M} = \bm{h} \nu\).
\end{enumerate}
We call \(\bm{h}\) the \emph{density function of \(\bm{M}\) with respect to the base \(\nu\)}.
\end{definition}

If \(\bm{M}\) is pure point,
we can take \(\nu\) to be the counting measure on \(\bR^d\)
and let \(\bm{h} : \bR^d \to \cM(\bR^d, \bC^{\ell\times\ell})\)
be given by \(\bm{h}(\{x\}) = \bm{M}(\{x\})\) for all \(x \in \bR^d\):
this is effectively the case already considered by \citet*{Baake_2019_origin}.
We will prove self-similar densities of \(\hat{\bm{H}}\) also exist
whenever \(\hat{\bm{H}}\) has a continuous part.

In order to study the diffraction measure \(\hat{\bm{H}}\),
we want to use a density \((\bm{h}, \nu)\) which is compatible with the self-similarity of \(\hat{\bm{H}}\).
Specifically, we want to find a density \((\bm{h}, \nu)\)
such that
the density function \(\bm{h}\) satisfies a renormalisation relation
analogous to the one satisfied by \(\hat{\bm{H}}\).
Recall that, for \(R > 0\),
\(B_R^\times = B_R \setminus \{0\}\)
is the punctured ball of radius \(R\) around \(0\).

\begin{definition}
\index[defs]{dilation invariant measure}
Let \(R > 0\).
\begin{enumerate}
    \item \index[defs]{dilation invariant measure}
    A positive measure \(\nu\) on \(B_R^\times\)
    is \emph{\(\lambda\)-dilation invariant} if \(D_\lambda \nu |_{B_R^\times} = \nu\).
    \item \index[defs]{self-similar density}
    A \emph{self-similar density of \(\hat{\bm{H}}\) on \(B_R^\times\)}
    is a density \((\bm{h}, \nu)\) of \(\hat{\bm{H}}|_{B_R^\times}\)
    such that the base \(\nu\) is \(\lambda\)-dilation invariant.
\end{enumerate}
\end{definition}

\begin{theorem}\label{thm:renormalisation-measure-existence}
For all \(R > 0\), \(\hat{\bm{H}}\)
has a self-similar density on \(B_R^\times\).
\end{theorem}
\begin{proof}
If \(\nu_1, \nu_2\) are two positive measures on a space \(X\),
we write \(\nu_1 \ll \nu_2\) if \(\nu_1\) is absolutely continuous with respect to \(\nu_2\),
i.e. \(\nu_1(S) = 0\) for all Borel sets \(S\) such that \(\nu_2(S) = 0\).
Our goal is to construct a \(\sigma\)-finite, \(\lambda\)-dilation invariant measure \(\nu\) on \(B_R\)
such that \(\norm{\hat{\bm{H}}}|_{B_R} \ll \nu\),
as then the Radon--Nikodym theorem implies that \(\hat{\bm{H}}|_{B_R}\)
has a density function \(\bm{h}\) with respect to the base \(\nu\),
so that \((\bm{h}, \nu)\) is a self-similar density of \(\hat{\bm{H}}\) on \(B_R^\times\).

By Theorem \ref{thm:diffraction-renormalisation-measure},
we have \(\norm{\hat{\bm{H}}} = \norm{\bm{A} (D_{\lambda^{-1}}\hat{\bm{H}}) \bm{A}^*} \leq \norm{\bm{A}}^2 \norm{D_{\lambda^{-1}}\hat{\bm{H}}}\),
where by \(\norm{\bm{A}}\) we mean the positive function \(\xi \mapsto \norm{\bm{A}(\xi)}\),
and by \(\norm{\hat{\bm{H}}}\) and \(\norm{D_{\lambda^{-1}}\hat{\bm{H}}}\) we mean the total variation measures of \(\hat{\bm{H}}\) and \(D_{\lambda^{-1}}\hat{\bm{H}}\).
Therefore \(\norm{\hat{\bm{H}}} \ll \norm{D_{\lambda^{-1}}\hat{\bm{H}}}\).
By iterating this process,
we deduce \(\norm{\hat{\bm{H}}} \ll \norm{D_{\lambda^{-n}}\hat{\bm{H}}}\) for all \(n \in \bN\).

Now let \(L = B_R \setminus B_{\lambda^{-1}R}\).
We define the positive measure \(\nu_0 := \norm{\hat{\bm{H}}}|_L\),
to be the total variation of \(\hat{\bm{H}}\) restricted to \(L\),
and we define a positive measure \(\nu\) by \[
    \nu = \sum_{n=0}^\infty D_{\lambda^{-n}} \nu_0
\]
As \(\norm{\hat{\bm{H}}}\) is a positive, \(\sigma\)-finite measure, so is \(\nu\).
By construction, it is clear that \(D_\lambda \nu|_{B_R^\times} = \nu\).
We want to show that \(\norm{\hat{\bm{H}}}|_{B_R} \ll \nu\),
as this implies the existence of a density function by the Radon--Nikodym theorem.

To see this, let \(S \subset B_R^\times\) be a Borel set such that \(\nu(S) = 0\)
and write \(S_n := S \cap D_{\lambda^{-n}}L\).
As \(\nu\) is positive, we have \(\nu(S_n) = 0\) for all \(n \in \bN_0\).
But we have \(\norm{\hat{\bm{H}}} \ll D_{\lambda^{-n}}\norm{\hat{\bm{H}}}\),
hence \(\nu(S_n) = D_{\lambda^{-n}}\nu_0(S_n) = 0\)
implies \(\norm{\hat{\bm{H}}}(S_n) = 0\).
Therefore, \(\norm{\hat{\bm{H}}}(S) = \sum_{n=0}^{\infty} \norm{\hat{\bm{H}}}(S_n) = 0\).
\end{proof}

Now we show that self-similar densities satisfy a renormalisation relation.

\begin{theorem}
Let \((\bm{h}, \nu)\) be a self-similar density of \(\hat{\bm{H}}\)
on \(B_R^\times\).
Then the density function \(\bm{h}\) satisfies \[
    \bm{h}(\xi) = \bm{A}(\xi)\bm{h}(\lambda \xi)\bm{A}(\xi)^*
\] for \(\nu\)-almost every \(\xi \in B_{\lambda^{-1}R}\).
\end{theorem}
\begin{proof}
Let \(f \in C_c(B_{\lambda^{-1}R})\).
Then \begin{align*}
D_{\lambda^{-1}}\hat{\bm{H}}(f) &= \int_{B_{R}} f(\lambda^{-1}\xi)\bm{h}(\xi) \,\nu(\mathrm{d}\xi)\\
    &= \int_{B_{\lambda^{-1}R}} f(\xi)\bm{h}(\lambda \xi) \,D_\lambda\nu(\mathrm{d}\xi) \\
    &= \int_{B_{\lambda^{-1}R}} f(\xi)\bm{h}(\lambda \xi) \,\nu(\mathrm{d}\xi)
\end{align*}

Hence the density of \(D_{\lambda^{-1}}\hat{\bm{H}}\) with respect to \(\nu\)
is given by \(\xi \mapsto \bm{h}(\lambda \xi)\).
We know \(\hat{\bm{H}}\) satisfies \(\hat{\bm{H}} = \bm{A}(\cdot) D_{\lambda^{-1}}\hat{\bm{H}} \bm{A}(\cdot)^*\),
so comparing the densities of both sides with respect to \(\nu\)
gives \[
    \bm{h}(\xi) = \bm{A}(\xi)\bm{h}(\lambda \xi)\bm{A}(\xi)^*
\]
for \(\nu\)-almost every \(\xi \in B_{\lambda^{-1}R}\),
which completes the proof.
\end{proof}

We want to use this to study
the decay of the diffraction measures of \(\cS\) around the origin.
This is given by the following lemma.
(We define \(\lesssim\), \(\gtrsim\), \(\lessapprox\) and \(\gtrapprox\)
in \Cref{sec:lyapunov-asymptotic-notation}).
\begin{theorem}[\citep*{thesis}]\label{lemma:uniformity-from-density}
Let \(\bm{\Phi}\) be a self-similar point process
with diffraction matrix \(\hat{\bm{H}}\),
whose substitution rule \(\cS\)
has scaling constant \(\lambda > 1\).

Let \(R > 0\)
and let \((\bm{h}, \nu)\) be a self-similar density
of \(\hat{\bm{H}}\)
on \(B_R^\times\).
For \(w \in \bC^\ell\),
we let \(\hat{\eta}_w = \langle \hat{\bm{H}}w, w\rangle\)
be the diffraction of \(\cS\) with weights \(w\).
\begin{enumerate}[(i)]
    \item Let \(\beta > 0\),
    and assume \[\langle \bm{h}(\lambda^{-N}\xi)w, w\rangle \lesssim  \norm{\bm{h}(\xi)}\lambda^{-\beta N} \]
    as \(N \to \infty\) uniformly for \(\nu\)-almost every \(\xi \in B_R\).
    Then we have \(\hat{\eta}_w(B_r) \lesssim r^{2\beta}\).
    \item Let \(\beta > 0\)
    and \(A \subset B_R\) be a Borel set
    such that \(\int_A \norm{\bm{h}(\xi)}\,\nu(\mathrm{d}\xi) > 0\),
    and assume \[\langle \bm{h}(\lambda^{-N}\xi)w, w\rangle \gtrsim \norm{\bm{h}(\xi)} \lambda^{-\beta N}\]
    as \(N \to \infty\) uniformly for \(\nu\)-almost every \(\xi \in A\).
    Then we have \(\hat{\eta}_w(B_r) \gtrsim r^{2\beta}\).
\end{enumerate}
The same statements hold
if we replace \(\lesssim\) with \(\lessapprox\)
and \(\gtrsim\) with \(\gtrapprox\).
\end{theorem}

In order to use \Cref{lemma:uniformity-from-density}, we want to understand \(\langle\bm{h}(\lambda^{-N}\xi)w,w\rangle\)
as \(N \to \infty\), using the renormalisation relation satisfied by \(\bm{h}\).
This motivates the following definition:
\begin{definition}

The \emph{Fourier cocycle} of the substitution rule \(\cS\)
is the matrix function given by \[
    \bm{A}^{(N)}(\xi) = \bm{A}(\lambda^{-N+1}\xi)\bm{A}(\lambda^{-N+2}\xi)\cdots \bm{A}(\xi)
\]
for \(\xi \in \bR^d\) and \(N \in \bN\).
\end{definition}

Then applying the renormalisation relation repeatedly yields the following corollary.
\begin{corollary}\label{corollary:renormalisation-density}
Let \((\bm{h}, \nu)\) be a self-similar density
of \(\hat{\bm{H}}\)
on \(B_R^\times\)
with respect to the substitution rule \(\cS\).
Then the density function \(\bm{h}\) satisfies \[
    \bm{h}(\lambda^{-N}\xi) = \bm{A}^{(N)}(\lambda^{-1}\xi)\bm{h}(\xi)\bm{A}^{(N)}(\lambda^{-1}\xi)^*
\]
for all \(\xi \in B_R^\times\), \(N \in \bN\).
\end{corollary}

\section{Linear cocycles around the origin}\label{chap:lyapunov}

In this section,
we prove Theorem \ref{thm:main-theorem-lyapunov},
which controls the behaviour of linear cocycles
such as the ones arising from the renormalisation relations of self-similar tilings.

\subsection{Asymptotic notation}
\label{sec:lyapunov-asymptotic-notation}
In order to state our results,
we will need some asymptotic notation,
which we define here.

\begin{definition}
Let \(X\) be a topological space, \(Y\) a set,
\(x_0 \in X\) and \(f, g : X \setminus \{x_0\} \times Y \to \bR_{\geq 0}\)
be two functions.
\begin{itemize}
    \item We write \(f(x, y) \lesssim g(x,y)\) or \(f(x,y) = O(g(x,y))\)
    \emph{as \(x \to x_0\) uniformly for \(y \in Y\)}
    if there exists a neighbourhood \(U\) of \(x_0\)
    and a constant \(C > 0\) such that
    \(f(x, y) \leq Cg(x,y)\) for all \(x \in U\), \(y \in Y\).
    We also write \(f(x,y) \gtrsim g(x,y)\)
    if \(g(x,y) \lesssim f(x,y)\),
    and \(\asymp\) if both \(\lesssim\) and \(\gtrsim\) hold.
    \item We write \(f(x,y) = o\!\left(g(x,y)\right)\)
    \emph{as \(x \to x_0\) uniformly for \(y \in Y\)}
    if for every \(\epsilon > 0\) there exists a neighbourhood \(U\) of \(x_0\)
    such that \(f(x,y) \le \epsilon\, g(x,y)\) for all \(x \in U\), \(y \in Y\).
\end{itemize}
\end{definition}

If \(X = \bN\), we define these asymptotic symbols as \(N \to \infty\)
in the usual way. We can also define these symbols
for functions \(f, g : X \setminus \{x_0\} \to \bR_{\geq 0}\) as \(x \to x_0\),
by letting \(Y\) be the set with one element:
in this way one recovers the usual meaning of \(O\), \(o\), \(\lesssim\) and \(\gtrsim\).

We also define a new asymptotic notation,
which is weaker than the usual asymptotic inequalities:
intuitively, it measures decay only \enquote{up to subexponential factors}.
\begin{definition}
Let \(Y\) be a set.
\begin{enumerate}
    \item Let \(f, g : \bN \times Y \to \bR_{\geq 0}\).
    We write \emph{\(f(N,y) \lessapprox g(N,y)\) as \(N \to \infty\) uniformly for \(y \in Y\)}
    if \(\limsup_{N \to \infty} \frac1N \log \sup_{y \in Y} f(N,y) \leq \limsup_{N \to \infty} \frac1N \log \sup_{y \in Y} g(N,y)\).
    \item Let \(f, g : \bR_{\geq 0} \times Y \to \bR_{\geq 0}\).
    We write \emph{\(f(r,y) \lessapprox g(r,y)\) as \(r \to 0\) uniformly for \(y \in Y\)}
    if \(\limsup_{r \to 0} \frac1{\log r} \log \sup_{y \in Y} f(r,y) \leq \limsup_{r \to 0} \frac1{\log r} \log \sup_{y \in Y}g(r,y)\).
\end{enumerate}
We define \(\gtrapprox\) analogously,
using \(\liminf\) instead of \(\limsup\),
and write \(\approx\) if both \(\lessapprox\) and \(\gtrapprox\) hold.
We can also define these notions without the second argument \(y\),
in which case we write \(f(N) \lessapprox g(N)\) as \(N \to \infty\).
Note that we set \(\log 0 = -\infty\).
\end{definition}

For example,
we have \(N \lessapprox 1\) as \(N \to \infty\),
as \(N\) grows slower than any exponential \(e^{\epsilon N}\) for \(\epsilon > 0\).
Other important special cases are:
\begin{itemize}
    \item For \(\alpha \in \bR\), \(f(r) \lessapprox r^\alpha\) as \(r \to 0\)
    if and only if, for all \(\epsilon > 0\), \(f(r) = O(r^{\alpha + \epsilon})\) as \(r \to 0\).
    \item \(f(r) \lessapprox 0\) as \(r \to 0\)
    if and only if, for all \(k \in \bR\), \(f(r) = O(r^k)\) as \(r \to 0\). 
\end{itemize}

\subsection{Linear cocycles around the origin}
\label{sec:lyapunov-exponents-converging-sequences}
Let \(\norm{\cdot}\) be any norm on \(\bC^{\ell}\),
and fix \(\lambda > 1\).
Let \(A : \bR^d \to \bC^{\ell \times \ell}\)
be a matrix-valued function:
unlike in the previous sections,
here we do not write it in bold,
in order to ease the notation.
For all \(\xi \in \bR^d\) and \(N \in \bN\)
we define the cocycle associated to \(A\) by
the product \(A^{(N)}(\xi) = A(\lambda^{-N+1}\xi) \cdots A(\xi) \).
In this section we prove Theorem \ref{thm:main-theorem-lyapunov},
which is a theorem about the asymptotic behaviour of linear cocycles.
We assume \(A\) is locally Lipschitz at \(0\),
i.e. \(A(\xi) = A(0) + O(\norm{\xi})\):
this holds in particular if \(A\) is smooth,
as is the case for the normalized Fourier matrix of a substitution rule.
In particular, \(A\) is \emph{locally Lipschitz at 0 with radius \(R\)}
if there exists \(L > 0\) such that \({\norm{A(\xi) - A(0)}}\leq L \norm{\xi}\) for all \(\xi \in B_R \setminus \{0\}\).

Let \(\chi_1 > \cdots > \chi_l\)
be the distinct values of \(\{\log \abs{\mu} \mid \mu \in \Spec A(0) \setminus \{0\}\}\)
and let \(\chi_{l+1} = -\infty\):
we call these the \emph{Lyapunov exponents} of \(A(0)\).
For \(\mu \in \Spec A(0)\), denote its generalized eigenspace
by \(E_\mu\),
and let \(E_j := \oplus \{E_\mu \mid \mu \in \Spec A(0), \log \abs{\mu} = \chi_j\}\)
be the space of generalized eigenvectors associated to \(\chi_j\)
for \(j \in [l+1]\).
Then there exist uniquely defined idempotent operators \(P_j : \bC^\ell \to \bC^\ell\)
for \(j \in [l+1]\)
such that \(\im P_j = E_j\) and \(\sum_{j=1}^{l+1}P_j = I\):
we call these the \emph{projections associated to the spaces \(E_j\).} 
Note that we are not assuming that the \(P_j\) are orthogonal projections.
Also define \(m_j := \dim E_j\).

We would like to know how \(A^{(N)}(\xi)x\) behaves as \(N \to \infty\).
If we restrict ourselves to the case of matrix powers,
we can obtain the following result
by elementary linear algebra.
\begin{example}
Assume \(B \in \bC^{\ell \times \ell}\) is unitarily diagonalizable,
let \(\chi_1, \cdots, \chi_l\)
be the finitely many values of \(\{\log \abs{\mu}\mid \mu \in \Spec B\}\),
and let \(E_1, \cdots, E_l\) be the corresponding sums of eigenspaces.
Then, if \(x \in E_j\), we have \begin{align*}
    \norm{B^N x} &= e^{\chi_j N} \norm{x} \\
    B^N x &\in E_j
\end{align*}
for all \(N \in \bN\).
\end{example}

Our goal is
to prove an analogous result for the linear cocycle \(A^{(N)}(\xi)\).
\begin{theorem}\label{thm:main-theorem-lyapunov}
There exists \(R > 0\)
such that the following holds:
for all \(\xi \in B_R\),
there exist projections \(P_j(\xi)\) for all \(j \in [l+1]\)
such that \(\sum_{j=1}^{l+1} P_j(\xi) = I\)
and the following asymptotic inequalities hold:
\begin{enumerate}
    \item For all \(j \in [l+1]\) we have
    \begin{equation}
             \norm{A^{(N)}(\xi)P_j(\xi)x} \approx e^{\chi_j N}\norm{P_j(\xi)x}\label{eq:lyapunov-norm}
    \end{equation}
    as \(N \to \infty\)
    uniformly for \(\xi \in B_R, x \in \bC^\ell\).
    \item For all \(j, k \in [l+1]\) we have
    \begin{align}
            {\norm{P_kA^{(N)}(\xi)P_j(\xi)x}} &\lesssim \lambda^{-N}\norm{A^{(N)}(\xi)P_j(\xi)x}
    & \mbox{if \(j > k\)} \label{eq:lyapunov-orth-forward}\\
    {\norm{P_kA^{(N)}(\xi)P_j(\xi)x}} &\lessapprox \max(e^{\chi_k - \chi_j}, \lambda^{-1})^N\norm{A^{(N)}(\xi)P_j(\xi)x}
        & \mbox{if \(j < k\)} \label{eq:lyapunov-orth-backward}
    \end{align}
    as \(N \to \infty\) uniformly for \(\xi \in B_R, x \in \bC^\ell\).
    \item For all \(j \in [l]\) such that \(m_j = 1\):
    \begin{equation}
            \norm{A^{(N)}(\xi)P_j(\xi)x} \gtrsim e^{\chi_j N} \norm{P_j(\xi)x} \label{eq:lyapunov-norm-mult-1}
    \end{equation}
        as \(N \to \infty\)
    uniformly for \(\xi \in B_R, x \in \bC^\ell\).
    \item If \(A(0)\) is diagonalizable,
    for all \(j \in [l]\):
        \begin{equation}
                \norm{A^{(N)}(\xi)P_j(\xi)x} \asymp e^{\chi_j N}\norm{P_j(\xi)x} \label{eq:lyapunov-norm-diagonalizable}
        \end{equation}
        as \(N \to \infty\)
        uniformly for \(\xi \in B_R, x \in \bC^\ell\).
\end{enumerate}
\end{theorem}

\begin{remark}
\citet*{Dubiner1991} proved a similar theorem
for the case of a sequence of matrices \(A_n\)
converging to a limit matrix \(A\),
without assumptions on the speed of convergence.
In this setting he was able to show Inequality \eqref{eq:lyapunov-norm}
and weaker versions of Inequalities \eqref{eq:lyapunov-orth-backward}
and \eqref{eq:lyapunov-orth-forward}.
However, for our application,
the stronger bounds on the speed of decay in Inequalities \eqref{eq:lyapunov-orth-backward} and \eqref{eq:lyapunov-orth-forward}
will be crucial,
as will be the uniformity for \(\xi\).
\end{remark}

\begin{remark}
In general, if \((A_n)_{n \in \bN}\) is an arbitrary sequence of matrices
and \(A^{(N)} = A_N \cdots A_1\) is its associated linear cocycle,
one asks for the possible values of the set
\[\{\lim_{N \to \infty} \frac1N \log \norm{A^{(N)}x} \mid x \in \bC^\ell\},\]
as well as the spaces where these values are attained.
These are called the \emph{Lyapunov exponents} and \emph{Lyapunov spaces} of the sequence,
respectively: therefore, our theorem characterizes the Lyapunov exponents and spaces
of the sequence \((A(\lambda^{-n}\xi))_{n \in \bN}\)
in terms of the Lyapunov exponents and spaces of the limit matrix \(A(0)\).

There exists a large literature on the theory of Lyapunov exponents and spaces,
particularly in the setting of cocycles arising from dynamical systems,
or their continuous analogues:
see \citep*{Barreira2017, Arnold1998}.
However, note that our case does not fit into the most common setting of Lyapunov theory:
the map \(x \mapsto \lambda^{-1}x\) is not measure preserving,
so our cocycle does not satisfy the conditions of Oseledets' Theorem.
This is why we do not make use of the existing literature on Lyapunov theory,
instead providing an elementary proof of Theorem \ref{thm:main-theorem-lyapunov}.
\end{remark}

\begin{remark}
If we equip \(\bC^\ell\) with an inner product
such that the projections \(P_j\) and \(P_j(\xi)\)
are orthogonal projections for all \(j \in [l+1]\),
we can interpret the estimates in the theorem via the formula \[
    \frac{\norm{P_kA^{(N)}(\xi)P_j(\xi)x}}{\norm{A^{(N)}(\xi)P_j(\xi)x}}
    = \cos \angle (E_k, A^{(N)}(\xi)P_j(\xi)x),
\]

where the angle between a subspace \(U\) and a vector \(v\)
is given by \({\angle(U,v) = \min_{u \in U} \angle(u,v)}\). Therefore, \Cref{eq:lyapunov-orth-backward,eq:lyapunov-orth-forward}
tell us that \(A^{(N)}(\xi)P_j(\xi)x\)
is \enquote{asymptotically orthogonal} to \(E_k\)
whenever \(k \neq j\):
the intuition is that \(A^{(N)}(\xi)P_j(\xi)x\) is close to being in \(E_j\),
which is orthogonal to the other eigenspaces.
In what follows we will mostly reason using norms and inner products instead of angles,
but it is still useful to keep this image in mind.
\end{remark}

\subsection{Exterior powers}
\label{sec:lyapunov-exterior-powers}

In order to prove Theorem \ref{thm:main-theorem-lyapunov},
we will use the formalism of exterior products
in order to represent linear subspaces.
Here we spell out some basic facts about exterior products which we need:
everything in this part is either known or elementary,
but we were not able to find a convenient reference for all of them.
The proofs can be found in \citep*[Chapter 5.2]{thesis}.

If \(V\) is a finite dimensional vector space over \(\bC\),
we denote by \(V^{\wedge q}\) the \(q\)-th exterior power of \(V\),
which can be defined as the dual of the space of alternating multilinear maps
from \(V^q\) to \(\bC\),
or as an appropriate quotient of the \(q\)-th tensor power \(V^{\otimes q}\):
see \citet*{kowalsky} for a more detailed definition.
For \(q_1, q_2 \in \bN\)
we denote the \emph{wedge product} by \(\wedge : V^{\wedge q_1} \times V^{\wedge q_2} \to V^{\wedge (q_1 + q_2)}\).
Then for all \(q \in \bN\),
the map \((v_1, \cdots, v_q) \mapsto v_1 \wedge \cdots \wedge v_q\)
is an alternating multilinear map.

For the rest of this section,
let \(\ell, q \in \bN\)
and fix an inner product on \(\bC^\ell\).
We consider the exterior power \((\bC^\ell)^{\wedge q}\),
which we denote as \(\bC^{\ell \wedge q}\) to ease the notation.
Then we can define a corresponding inner product on \((\bC^\ell)^{\wedge q}\)
with the following property:
if \(\alpha = v_1 \wedge \cdots \wedge v_q, \beta = w_1 \wedge \cdots \wedge w_q\),
\[
    \langle \alpha, \beta \rangle = \det J(\alpha, \beta)
\]
where \(J(\alpha, \beta) \in \bC^{\ell \times \ell}\)
is the matrix with entries \(J_{i,k} = \langle v_i, w_k \rangle\).
This inner product makes the wedge product continuous.

If \(\alpha \in \bC^{\ell \wedge q}\)
can be written as \(\alpha = a_1 \wedge \cdots \wedge a_q\)
for linearly independent vectors \(a_1, \cdots, a_q \in \bC^\ell\),
we say \(\alpha\) is a \emph{\(q\)-blade} or \emph{decomposable \(q\)-vector},
and let its \emph{span} be defined by \(\linspan \alpha = \linspan (a_1, \cdots, a_q)\),
where the \(a_j\) are any \(q\) vectors such that \(\alpha = a_1 \wedge \cdots \wedge a_q\).
(This quantity does not depend on the choice of the \(a_j\).)
We also let \(\pi_\alpha\) be the orthogonal projection onto the subspace \(\linspan \alpha \subset \bC^\ell\).

For \(q_1, q_2 \in \bN\),
\(\alpha \in \bC^{\ell \wedge q_1}\),
and \(\beta \in \bC^{\ell \wedge q_2}\),
we write \(\alpha \subset \beta\) if there exists some \(\gamma\)
such that \(\beta = \alpha \wedge \gamma\).
If \(\alpha, \beta\) are blades,
this is the case if and only if \(\linspan \alpha \subset \linspan \beta\).

We need the following properties of blades.
\begin{lemma}\label{lemma:compose-blade}
Let \(q, q' \in \bN\),
\(\alpha, \beta\) be \(q\)-blades, \(\alpha', \beta'\) be \(q'\)-blades.
\begin{enumerate}[(i)]
    \item Assume \(\linspan \alpha \perp \linspan \alpha'\).
    Then we have \(\langle \alpha \wedge \alpha', \alpha \wedge \beta'\rangle = \norm{\alpha}^2\langle \alpha', \beta'\rangle\).
    \item Assume \(\linspan \alpha \perp \linspan \alpha'\).
    Then we have \(\norm{\alpha \wedge \alpha'} = \norm{\alpha}\norm{\alpha'}\)
    \item For all \(x \in \bC^\ell\), we have \(\norm{\alpha \wedge x} = \norm{\alpha} \norm{(I-\pi_{\alpha})x}\).
\end{enumerate}
\end{lemma}
\begin{lemma}\label{lemma:blades-closed}
Let \(B_q(\bC^\ell) \subset \bC^{\ell \wedge q}\)
be the set of \(q\)-blades in \(\bC^{\ell \wedge q}\).
Then its closure in the norm topology of \(\bC^{\ell \wedge q}\)
is \(B_q(\bC^\ell) \cup \{0\}\).
\end{lemma}

\begin{lemma}\label{lemma:projection-blade-bound}
For all \(q\)-blades \(\alpha, \beta\)
we have
\[\norm{\pi_\alpha - \pi_\beta} \leq 2 \frac{\norm{\alpha - \beta}}{\min(\norm{\alpha},\norm{\beta})}
\]
\end{lemma}

Now let \(A \in \bC^{\ell \times \ell}\).
We define the \emph{\(q\)-th exterior power of \(A\)}
as the linear map \(A^{\wedge q} : \bC^{\ell \wedge q} \to \bC^{\ell \wedge q}\)
by \[
    A^{\wedge q}(v_1 \wedge \cdots \wedge v_q) = A v_1 \wedge \cdots \wedge A v_q
\]
for all \(v_1, \cdots, v_q \in \bC^\ell\).
Then the generalized eigenvalues of \(A^{\wedge q}\)
are given by the generalized eigenvalues of \(A\).
\begin{lemma}
\label{lemma:exterior-power-generalized-eigenvalues}

Let \(A \in \bC^{\ell \times \ell}\)
and let \(B = (v_1, \cdots, v_{\ell})\)
be an ordered basis of \(\bC^\ell\)
such that \(A\) is in upper Jordan normal form
with respect to this basis.
Let \(\mu_1, \cdots, \mu_\ell\) be the corresponding diagonal entries of \(A\),
which are the generalized eigenvalues of \(A\)
repeated according to their multiplicity.

Let \[B' = \{v_{j_1} \wedge \cdots \wedge v_{j_q}
\mid j_1, \cdots, j_q \in [\ell], j_1 < \cdots < j_q\}
\]
be the corresponding basis of \(\bC^{\ell \wedge q}\),
ordered lexicographically on the indices.
Let \((a_{JK})_{J, K}\)
be the coefficients of the transformation matrix of \(A^{\wedge q}\)
with respect to the basis \(B'\),
indexed by ordered \(q\)-tuples of indices \(J, K\).
Then \(A^{\wedge q}\) has the following properties:
\begin{enumerate}
    \item For \(J = (j_1, \cdots, j_q)\),
    the corresponding diagonal entry of \(A^{\wedge q}\) is given by \(a_{JJ} = \mu_{j_1} \cdots \mu_{j_q}\).
    \item For \(J = (j_1, \cdots, j_q)\), \(K = (k_1, \cdots, k_q)\)
    with \(J < K\), we have \begin{itemize}
        \item \(a_{KJ} = 0\)
        \item \(a_{JK} \neq 0 \implies a_{JJ} = a_{KK}\).
    \end{itemize}
    \item \(A^{\wedge q}\) is upper triangular with respect to the basis \(B'\),
    and its generalized eigenvalues with multiplicity
    are given by \(\{\mu_{j_1} \cdots \mu_{j_q} \mid 1 \leq j_1 < \cdots < j_q \leq \ell\}\).
\end{enumerate}
\end{lemma}

\begin{lemma}
\label{lemma:exterior-power-continuity}
For all \(A, B \in \bC^{\ell \times \ell}\),
the inequality \[
    \norm{A^{\wedge q} - B^{\wedge q}} \leq q \max(\norm{A}, \norm{B})^{q-1} \norm{A - B}
\]
holds.
\end{lemma}

\subsection{Proof of Theorem \ref{thm:main-theorem-lyapunov}}
\label{sec:lyapunov-proof}

Now we set out to prove Theorem \ref{thm:main-theorem-lyapunov}.
The proof is morally analogous to the QR algorithm
for computing eigenvalues,
which has been generalized to compute Lyapunov exponents
\citep*{Ginelli_2013, Noethen_2019}.

\subsection{The power method}
For this section, equip \(\bC^\ell\) with an arbitrary inner product
and its corresponding norm:
in the next section, we will make a specific choice of inner product.

\begin{lemma}\label{prod-lq-sum}
Let \((a_m)_{m \in \bN}\) be a sequence of positive real numbers.
Then, for all \(M \in \bN\),
the following inequality holds\[
    \prod_{m=1}^M a_m \leq \exp\left(\sum_{m=1}^M (a_m - 1)\right)
\]
\end{lemma}
\begin{proof}
Take the logarithm of both sides of the inequality
and use \(\log x = \log(1 + (x-1)) \leq x-1\).
\end{proof}

\begin{lemma}\label{lemma:cocycle-norm-lq}
Let \(A : \bR^d \to \bC^{\ell \times \ell}\) be a matrix function
which is locally Lipschitz at \(0\) with radius \(R\).
Then the following inequalities hold:
\begin{enumerate}
    \item For all \(k \in \bN\) such that \(A(0)^k \neq 0\),
the matrix cocycle satisfies \[
    \norm{A^{(N)}(\xi)} \lesssim \norm{A(0)^k}^{\frac{N}{k}}
\]
as \(N \to \infty\) uniformly for \(\xi \in B_R\).
\item The matrix cocycle satisfies \[
    \norm{A^{(N)}(\xi)} \lessapprox \rho(A(0))^N
\] as \(N \to \infty\) uniformly for \(\xi \in B_R\),
where \(\rho(A(0)) = \max \{\abs{\mu} \mid \mu \in \Spec A(0)\}\)
is the spectral radius of \(A(0)\).
\end{enumerate}
\end{lemma}
\begin{proof}
We prove both parts separately:
\begin{enumerate}
    \item For the first part, let \(N \in \bN\):
    we can always write it as \(N = kM + r\) for some \(M \in \bN\), \(r \in \{0, \cdots, k-1\}\).
    Then taking the limit as \(N \to \infty\) is the same as taking the limit as \(M \to \infty\)
    for all \(r \in \{0, \cdots, k-1\}\).

Recall that, for \(N \in \bN\),
the linear cocycle is given by \(A^{(N)}(\xi) = A(\lambda^{-N+1}\xi) \cdots A(\xi)\).
Then we can decompose it as \[
    A^{(N)}(\xi) = A^{(r)}(\lambda^{-kM}\xi) \prod_{m=0}^{M-1} A^{(k)}(\lambda^{-km} \xi)
\]
(where the iterated product is taken from right to left)
because
\begin{align*}
A^{(r)}(\lambda^{-kM}\xi) &= A(\lambda^{-(kM+r)+1}\xi)\cdots A(\lambda^{-kM}\xi)\\
A^{(k)}(\lambda^{-km}\xi) &= A(\lambda^{-k(m+1)+1}\xi)\cdots A(\lambda^{-km}\xi)
\end{align*}
Therefore, using the submultiplicativity of the matrix norm,
we obtain \[
    \frac{\norm{A^{(N)}(\xi)}}{\norm{A(0)^k}^M} 
    \leq \norm{A^{(r)}(\lambda^{-kM}\xi)}
    \prod_{m=0}^{M-1} \frac{\norm{A^{(k)}(\lambda^{-km} \xi)}}{\norm{A(0)^k}} 
\]

Using \(A(\xi) = A(0) + O(\norm{\xi})\) as \(\xi \to 0\)
and \Cref{lemma:exterior-power-continuity},
we can bound the two factors as follows.
\begin{itemize}
    \item We have \(\sum_{m=0}^{M-1} \left(\frac{\norm{A^{(k)}(\lambda^{-km} \xi)}}{\norm{A(0)^k}} - 1 \right)
\lesssim \norm{\xi}\) as \(M \to \infty\) uniformly for \({\xi \in B_R}\),
because we can bound the sum by a geometric series.
Therefore, by \Cref{prod-lq-sum},
\(\prod_{m=0}^{M-1} \frac{\norm{A^{(k)}(\lambda^{-km} \xi)}}{\norm{A(0)^k}} \lesssim e^{\norm{\xi}} \lesssim 1\)
as \(M \to \infty\) uniformly for \(\xi \in B_R\).
\item 
We have \({\norm{A^{(r)}(\lambda^{-kM}\xi)} \lesssim 1}\)
as \(M \to \infty\) uniformly for \(r \in \{0, \cdots, k-1\}\), \(\xi \in B_R\).
\end{itemize}

Putting both of these facts together, we have \[
    \frac{\norm{A^{(N)}(\xi)}}{\norm{A(0)^k}^M} 
    \leq \norm{A^{(r)}(\lambda^{-kM}\xi)}
    \prod_{m=0}^{M-1} \frac{\norm{A^{(k)}(\lambda^{-km} \xi)}}{\norm{A(0)^k}} 
    \lesssim 1
\]
as \(N \to \infty\) uniformly for \(\xi \in B_R\).

Finally, we have \(M = \frac{N}{k} - \frac{r}{k}\)
which means \[
    \norm{A^{(N)}(\xi)} \lesssim \norm{A(0)^k}^M = \norm{A(0)^k}^{\frac{N}{k}}\norm{A(0)^k}^{-\frac{r}{k}}
    \lesssim \norm{A(0)^k}^{\frac{N}{k}}
\]
as \(N \to \infty\) uniformly for \(\xi \in B_R\).
This concludes the proof of \textit{(i)}.

    \item For the proof of \textit{(ii)},
    we distinguish two cases:
    \begin{itemize}
        \item Assume \(A(0)^k = 0\) for some \(k \in \bN\):
        then \(\rho(A(0)) = 0\).
        This means we want to prove \(\norm{A^{(N)}(\xi)}\) decays faster than any exponential function.
        As in the first part, for all \(N \in \bN\),
        we can write \(N = kM + r\) with \(M \in \bN\), \(r \in \{0, \cdots, k-1\}\).
        Then we have \begin{align*}
            \frac1N \log \norm{A^{(N)}(\xi)} &= \frac1N \log \norm{A^{(r)}(\lambda^{-kM}\xi)} + \frac1N \sum_{m=0}^{M-1} \log \norm{A^{(k)}(\lambda^{-km} \xi)} \\
            &\lesssim -\frac1N \sum_{m=0}^{M-1} m \lesssim -N
        \end{align*}
        This implies that \(\frac1N \log \norm{A^{(N)}(\xi)} \to -\infty\)
        as \(N \to \infty\), uniformly for \(\xi \in B_R\),
        therefore \(\norm{A^{(N)}(\xi)} \lessapprox 0 = \rho(A(0))\) as \(N \to \infty\) uniformly for \(\xi \in B_R\).
        \item Now assume \(A(0)\) is not nilpotent,
        so we can apply part \textit{(i)} for every \(k \in \bN\).
        Define a sequence \((a_n)_{n \in \bN}\) by \[
            a_N := \sup_{\xi \in B_R} \log \norm{A^{(N)}(\xi)}
        \] By the submultiplicativity of the matrix norm,
    this sequence is subadditive,
    so Fekete's Lemma tells us that
    \(\lim_{N \to \infty} \frac{a_N}{N} = \inf_{N \in \bN} \frac{a_N}{N}\).
    Let \(\chi := \lim_{N \to \infty} \frac{a_N}{N} = \lim_{N \to \infty} \frac1N \log \sup_{\xi \in B_R} \norm{A^{(N)}(\xi)}\).
    Then, by part \textit{(i)},
    for all \(k \in \bN\)
    we have \[
        \chi = \lim_{M \to \infty} \frac{a_{kM}}{kM} \leq \frac1k \log \norm{A(0)^k}
    \]
    By taking the limit as \(k \to \infty\),
    we obtain \(\chi \leq \lim_{k \to \infty} \frac1{k}\log \norm{A(0)^k} = \log \rho(A(0))\),
    which concludes the proof.
    \end{itemize}
\end{enumerate}
\end{proof}

\begin{lemma}\label{lemma:top-backwards-vector}
Let \(A : \bR^d \to \bC^{\ell \times \ell}\) be a matrix function
which is locally Lipschitz at \(0\) with radius \(R\),
and assume there exists a right eigenvector \(v\) of \(A(0)^*\)
whose eigenvalue \(\overline{\mu}\)
satisfies \(\abs{\overline{\mu}}^k = \norm{A(0)^k} > 0\)
for some \(k \in \bN\). We denote the conjugate of \(\overline{\mu}\) by \(\mu\).

Define \(w^{(N)}(\xi) := \overline{\mu}^{-N}A^{(N)*}(\xi)v\).
Then:
\begin{enumerate}[(i)]
    \item The limit \(w(\xi) := \lim_{N \to \infty} w^{(N)}(\xi)\) exists
        for all \(\xi \in B_R\).
    \item \(\norm{w(\xi) - w^{(N)}(\xi)} = O(\lambda^{-N}\norm{\xi})\)
    \item \(\norm{w(\xi) - v} = O(\norm{\xi})\)
    \item \(\norm{w(\xi)} = 1 + O(\norm{\xi})\)
    \item \(\langle A^{(N)}(\xi)w(\xi), v\rangle = \mu^N (1+O(\norm{\xi}))\)
    \item \(\norm{A^{(N)}(\xi)w(\xi)} \asymp \abs{\mu}^N\)
\end{enumerate}
as \(\xi \to 0\), \(N \to \infty\).
\end{lemma}
\begin{proof}
By Lemma \ref{lemma:cocycle-norm-lq},
we have \(\norm{A^{(N)}(\xi)} \lesssim \norm{A(0)^k}^{\frac{N}{k}} = \abs{\mu}^N\)
as \(N \to \infty\):
therefore \(\norm{\overline{\mu}^{-N-1}A^{(N)*}(\xi)} \lesssim 1\)
as \(\xi \to 0, N \to \infty\).

Then the differences \(w^{(N+1)}(\xi) - w^{(N)}(\xi)\)
are bounded by
\begin{align*}
    \norm{w^{(N+1)}(\xi) - w^{(N)}(\xi)} &= \norm{\overline{\mu}^{-N-1}A^{(N+1)*}(\xi)v - \overline{\mu}^{-N}A^{(N)*}(\xi)v} \\
    &\leq \norm{\overline{\mu}^{-N-1}A^{(N)*}(\xi)}\norm{A(\lambda^{-N}\xi)^*v - \overline{\mu} v} \\
    &\leq \norm{\overline{\mu}^{-N-1}A^{(N)*}(\xi)}\norm{A(\lambda^{-N}\xi)^* - A(0)^*}\norm{v} \\
    &= O(\lambda^{-N}\norm{\xi})
\end{align*}
as \(\xi \to 0, N \to \infty\).

We can write \(w^{(N)}(\xi) = v + \sum_{n=1}^{N} \left(w^{(n)}(\xi) - w^{(n-1)}(\xi)\right)\):
as \(N \to \infty\), the latter sum is bounded by a telescoping series with exponent \(\lambda\),
therefore the limit \(w(\xi)\) exists and satisfies \textit{(i)} and \textit{(iii)}.
Furthermore, we have \[
    w(\xi) - w^{(N)}(\xi) = \sum_{n=N+1}^\infty w^{(n)}(\xi) - w^{(n-1)}(\xi) = O(\lambda^{-N}\norm{\xi})
\]
as \(N \to \infty, \xi \to 0\), where again we are bounding the series by a geometric series.
Statement \textit{(iv)} follows directly from \textit{(iii)}.

For \textit{(v)}, we have 
\begin{align*}
    \langle A^{(N)}(\xi)w(\xi), v\rangle
    &= \langle w(\xi), \overline{\mu}^N \overline{\mu}^{-N} A^{(N)}(\xi)^*v\rangle \\
    &= \mu^N \langle w(\xi), w^{(N)}(\xi)\rangle \\
    &= \mu^N (\norm{w(\xi)}^2 + \langle w(\xi), w^{(N)}(\xi) - w(\xi)\rangle) \\
    &= \mu^N(1 + O(\norm{\xi}))
\end{align*}
as \(\xi \to 0, N \to \infty\), using \textit{(ii)} and \textit{(iv)}.

Now we only have to prove \textit{(vi)}.
From Lemma \ref{lemma:cocycle-norm-lq} we obtain
\(\norm{A^{(N)}(\xi)} \lesssim \abs{\mu}^N\),
which is the upper bound. 
For the lower bound, we use \textit{(v)}:
\begin{align*}
    \norm{A^{(N)}(\xi)w(\xi)} &\geq \abs{\langle A^{(N)}(\xi)w(\xi), v\rangle} \\
    &= \abs{\mu}^N(1 + O(\norm{\xi})) \qedhere
\end{align*}
\end{proof}

\begin{example}
If \(A(\xi)\) is the normalized Fourier matrix of a primitive substitution rule,
the vector given by \(v_j = m_{\bR^d}(\tau_j)\)
is an eigenvector of \(A(0)^*\) with eigenvalue \(1 = \norm{A(0)}\).
Then one can show that \(w(\xi)\) is given by
\(w(\xi)_j = \widehat{\bone_{\tau_j}}(\xi)\) (up to a constant factor),
as this is the vector function that satisfies \(w(\xi) = A(\xi)^*w(\lambda^{-1}\xi)\).
Then \(w(\xi)\) can be understood as the cocycle analogue of an eigenvector,
and the proof of Lemma \ref{lemma:top-backwards-vector}
is analogous to the power method for finding the top eigenvector.
\end{example}

\subsection{Defining the subspaces}
For the rest of this section,
we fix a matrix-valued function \(A\)
which is locally Lipschitz at \(0\).
Let \(\chi_j\), \(E_j\), \(P_j\) be as in \Cref{sec:lyapunov-exponents-converging-sequences}
and define \(m_j := \dim E_j\).

We will want to pick a basis of \(\bC^\ell\)
of the following form.
\begin{definition}
\index[defs]{descending Jordan basis}
A \emph{descending Jordan basis} for a matrix \(A \in \bC^{\ell \times \ell}\)
is a basis \(v_1, \cdots, v_\ell\) with the following properties:
\begin{itemize}
    \item \(A\) is in lower Jordan normal form with respect to this basis.
    \item The diagonal entries \(\mu_j\) with respect to the basis
    satisfy \(\abs{\mu_1} \geq \abs{\mu_2} \geq \cdots \geq \abs{\mu_\ell}\).
\end{itemize}
\end{definition}

Every square matrix over the complex numbers
has a Jordan normal form,
so in particular it admits a descending Jordan basis.
Furthermore, given a descending Jordan basis,
we can choose an appropriate inner product on \(\bC^\ell\)
which makes this basis orthonormal.
\textbf{For the rest of this section,
we fix a descending Jordan basis \(v_1, \cdots, v_\ell\)  for \(A(0)\)
and an inner product on \(\bC^\ell\)
such that the basis is orthonormal with respect to the inner product.}
Let \(\mu_1, \cdots, \mu_\ell\)
be corresponding diagonal entries of \(A(0)\).

By the definition of a descending Jordan basis,
the first \(m_1\) vectors span \(E_1\),
the next \(m_2\) span \(E_2\), and so on.
For \(j \in \{0, \cdots, l+1\}\),
define \(m_{(j)} = m_1 + \cdots + m_j\) as the sum of the first \(j\) multiplicities
for \(j = 0, \cdots, l+1\), setting \(m_{(0)} := 0\).
As the basis is orthonormal,
\(A(0)^*\) is in upper Jordan normal form with respect to this basis,
and the diagonal entries are given by \(\overline{\mu}_j\).

We will characterize the Lyapunov spaces of the cocycle
by using exterior algebra.
Recall that,
for \(A \in \bC^{\ell \times \ell}\),
\(A^{\wedge q}\) is the \emph{\(q\)-th exterior power} of \(A\),
which is the linear operator on \(\bC^{\ell \wedge q}\)
which satisfies \[
    A^{\wedge q}(w_1 \wedge \cdots \wedge w_q) = Aw_1 \wedge \cdots \wedge Aw_q
\]
for all \(w_1, \cdots, w_q \in \bC^\ell\).

In our case,
for all \(q \in [\ell]\),
we can consider the matrix function \(A^{\wedge q}\)
defined by \(A^{\wedge q}(\xi) = A(\xi)^{\wedge q}\),
and its associated cocycle
\(A^{(N)\wedge q}(\xi) = A(\lambda^{-N+1}\xi)^{\wedge q} \cdots A(\xi)^{\wedge q}\).

Then we can define:
\begin{definition}
Let \(j \in [\ell]\).
We define \begin{align*}
    \alpha_{j} &= v_1 \wedge \cdots \wedge v_{j} \\
    \mu_{(j)} &= \mu_1 \mu_2 \cdots \mu_{j} \\
    \beta_{j}^{(N)}(\xi) &= \overline{\mu_{(j)}}^{-N} A^{(N)*\wedge j}(\xi) \alpha_{j} \\
    \beta_{j}(\xi) &= \lim_{N \to \infty} \beta_{j}^{(N)}(\xi)
\end{align*}
whenever the latter limit exists.
Here, by \(A^{(N)*\wedge {j}}\)
we mean the product \((A^{(N)\wedge j}(\xi))^* = A^{\wedge j}(\xi)^* \cdots A^{\wedge j}(\lambda^{-N+1}\xi)^*\).
\end{definition}

We want to use Lemma \ref{lemma:top-backwards-vector}
to prove that,
for appropriate choices of \(q\),
the limits \(\beta_q(\xi)\) exist.
We do this by choosing \(q = m_{(j)}\)
for \(j = 1, \cdots, l+1\):
then \(\alpha_q\) is an eigenvector satisfying the assumptions of Lemma \ref{lemma:top-backwards-vector}.
\begin{lemma}\label{lemma:power-assumption}
Let \(j \in [l]\)
and \(q := m_{(j)}\).
Then the following properties hold:
\begin{enumerate}
    \item The matrix function \(A^{\wedge q}\) is locally Lipschitz at \(0\):
    that is, we have the bound \[{\norm{A(\xi)^{\wedge q} - A(0)^{\wedge q}} = O(\norm{\xi})}\]
    as \(\xi \to 0\).
    \item \(\alpha_{q}\) is an eigenvector of \(A(0)^{\wedge q*}\) with eigenvalue \(\overline{\mu}_{(q)}\).
    \item There exists \(k \in \bN\)
    such that \({\norm{(A(0)^{\wedge q})^k}} = \abs{\mu_{(q)}}^k\).
\end{enumerate}
\end{lemma}
\begin{proof}
Part \textit{(i)} follows from Lemma \ref{lemma:exterior-power-continuity},
using the fact that \(A\) itself is locally Lipschitz at \(0\).

For the rest, note that the structure of \(A(0)^{\wedge q*}\)
is described by \Cref{lemma:exterior-power-generalized-eigenvalues},
as the vectors \(v_1, \cdots, v_q\)
put \(A(0)^*\) in upper Jordan normal form by assumption:
this means the \(A(0)^{\wedge q*}\) is in upper triangular form,
where the diagonal entries are products of the eigenvalues of \(A(0)^*\) with multiplicity,
and the off-diagonal entries are associated to vectors in the same Jordan block of \(A(0)^*\).

By the way we have chosen \(q\), the vectors \(v_1, \cdots, v_q\)
are the generalized eigenvectors of \(A(0)^*\)
associated to the exponents \(\chi_1, \cdots, \chi_j\),
so that \(\abs{\mu_t} < e^{\chi_j}\) for all \(t > q\).
In particular, if \(s \leq q\) and \(t > q\), \(v_s\) and \(v_t\)
correspond to different Jordan blocks in \(A(0)^*\).
Then, using \Cref{lemma:exterior-power-generalized-eigenvalues},
the transformation matrix of \(A(0)^{\wedge q*}\)
can be written as \[
    \begin{pmatrix}
        \mu_{(q)} & \bm{0}^\top \\
        \bm{0} & M
    \end{pmatrix}
\]
where \(\bm{0}\) denotes the column vector of zeros of appropriate size,
and \(M\) is a matrix with spectral radius strictly less than \(\abs{\mu_{(q)}}\).

Therefore,
for any sufficiently large \(k \in \bN\),
we have \(\norm{M^k} < \abs{\mu_{(q)}}^k\)
and therefore \(\norm{(A(0)^{\wedge q*})^k} = \abs{\mu_{(q)}}^k\),
which concludes the proof of \textit{(iii)}.
\end{proof}

\begin{theorem}\label{thm:lyapunov-blades}
Let \(j \in [l]\),
set \(q := m_{(j)}\).
\begin{enumerate}[(i)]
    \item There exists some \(R > 0\)
    such that the limit \(\beta_{q}(\xi) := \lim_{N \to \infty}\beta_{q}^{(N)}(\xi)\)
    exists and is a \(q\)-blade or zero for all \(\xi \in B_R\).
    \item \(\norm{\beta_{q}^{(N)}(\xi) - \beta_{q}(\xi)} = O(\lambda^{-N}\norm{\xi})\)
    as \(\xi \to 0\).
    \item \(\norm{\beta_{q}(\xi) - \alpha_{q}} = O(\norm{\xi})\)
    as \(\xi \to 0\).
    \item \(\norm{\beta_{q}(\xi)} = 1 + O(\norm{\xi})\)
    as \(N \to \infty\), \(\xi \to 0\).
    \item \(\langle A^{(N)\wedge q}(\xi)\beta_{q}(\xi), \alpha_{q}\rangle = {\mu_{(q)}}^N(1 + O(\norm{\xi}))\)
    as \(N \to \infty\), \(\xi \to 0\).
    \item \(\norm{A^{(N)\wedge {q}}(\xi)\beta_{q}(\xi)} \asymp \abs{\mu_{(q)}}^N\)
    as \(N \to \infty\), \(\xi \to 0\).
\end{enumerate}
\end{theorem}
\begin{proof}
By Lemma \ref{lemma:power-assumption},
the matrix function \(A^{\wedge q}\)
satisfies all the assumptions of Lemma \ref{lemma:top-backwards-vector},
with respect to the vector \(\alpha_{q}\):
from this, almost every statement follows directly by taking \(v = \alpha_{q}\).
The only thing we have not proved yet is that \(\beta_{q}(\xi)\) is a \(q\)-blade or zero,
but this follows from the fact that \(\beta_{q}^{(N)}(\xi)\) is a blade for all \(N\)
and \Cref{lemma:blades-closed}.
\end{proof}

\begin{definition}
\label{def:lyapunov-projections}
For \(j \in [l+1]\), define the subspaces
\[
    Y_j(\xi) := \left(\linspan \beta_{m_{(j)}}(\xi)\right)^{\perp},
\]
and set \(Y_0(\xi) = \bC^\ell\).
Define the spaces
\[
    E_j(\xi) := Y_{j-1}(\xi) \cap Y_{j}(\xi)^\perp
\]
for \(j \in [l+1]\), and let \(P_j(\xi)\) be the orthogonal projection onto \(E_j(\xi)\).
We also set \(P_{(k)} := P_1 + \cdots + P_k\)
and \(P_{(k)}(\xi) := P_1(\xi) + \cdots + P_k(\xi)\).
\end{definition}

\subsection{Proving the inequalities}
Our goal will be to show that the projections \(P_j(\xi)\)
we just defined
satisfy the inequalities from Theorem \ref{thm:main-theorem-lyapunov}.
In particular,
this will mean that the subspaces \(Y_j(\xi)\) are the \emph{Lyapunov subspaces} of the linear cocycle
of the sequence \(A(\lambda^{-N+1}\xi)\) (see \citep*{Barreira2017} for a definition).

First, we make the following crucial observation:
\begin{lemma}
Let \(j \in [l+1]\).
Then \[
    A(\xi) Y_j(\xi) \subseteq Y_{j}(\lambda^{-1}\xi)
\]
\end{lemma}
\begin{proof}
Set \(q := m_{(j)}\)
and let \(x \in Y_j(\xi)\).
By definition, this means that for all \(u \in \linspan \beta_q(\xi)\),
we have \(\langle x, u \rangle = 0\).
Furthermore, by definition of the \(\beta_q\),
we know that \(A(\xi)^* \beta_q(\lambda^{-1}\xi)\) and \(\beta_q(\xi)\) are linearly dependent,
so \(\linspan A(\xi)^* \beta_q(\lambda^{-1}\xi) = \linspan \beta_q(\xi)\),
Now let \(u' \in \linspan \beta_q(\lambda^{-1}\xi)\):
then \(A(\xi)^* u' \in \linspan A(\xi)^* \beta_q(\lambda^{-1}\xi) = \linspan \beta_q(\xi)\),
which implies \(\langle A(\xi)x, u'\rangle = \langle x, A(\xi)^* u' \rangle = 0\):
and thus \(A(\xi)x \in \linspan \beta_q(\xi)^\perp\),
which is what we wanted to show.
\end{proof}

Using this, we are ready to prove the bounds in Theorem \ref{thm:main-theorem-lyapunov}.

\begin{lemma}\label{lemma:orthogonality-before}
Let \(j, k \in [l+1]\),
such that \(j < k\).
Then \[
\norm{P_j A^{(N)}(\xi)P_k(\xi)x} \lesssim \lambda^{-N} \norm{\xi} \norm{\bm{A}^{(N)}(\xi)P_k(\xi)x}
\]
as \(N \to \infty, \xi  \to 0\)
uniformly for \(x \in \bC^\ell\).
\end{lemma}
\begin{proof}
Let \(y = P_k(\xi)x\).
Without loss of generality we can assume that \(y \neq 0\),
as otherwise we have \(0\) on both sides of the inequality.
Also let \(q := m_{(j)}\).

As \(j < k\), we have \(y \in Y_j(\xi)\).
By definition of the Lyapunov subspaces,
we have \(A^{(N)}(\xi)Y_{j}(\xi) \subset Y_{j}(\lambda^{-N}\xi)\),
so \(A^{(N)}(\xi)y \in Y_{j}(\lambda^{-N}\xi)\).
We can rewrite this using orthogonal projections: \[
    A^{(N)}(\xi)y = (I - P_{(j)}(\lambda^{-N}\xi))A^{(N)}(\xi)y
\]

By \Cref{thm:lyapunov-blades},
we have \(\norm{P_m - P_m(\xi)} = O(\norm{\xi})\)
for all \(m \in [l]\),
which implies \(\norm{P_{(j)} - P_{(j)}(\xi)} = O(\norm{\xi})\),
therefore \(\norm{P_j(I - P_{(j)}(\xi))} = O(\norm{\xi})\)
as \(\xi \to 0\).
Therefore \begin{align*}
    \frac{\norm{P_j A^{(N)}(\xi)y}}{\norm{A^{(N)}(\xi)y}}
    &= \frac{\norm{P_j (I - P_{(j)}(\lambda^{-N}\xi))A^{(N)}(\xi)y}}{\norm{A^{(N)}(\xi)y}} \\
    &\leq \norm{P_j (I - P_{(j)}(\lambda^{-N}\xi))} \\
    &= O(\lambda^{-N}\norm{\xi})
\end{align*}
as \(N \to \infty, \xi\) uniformly on \(x \in \bC^\ell\).
\end{proof}

\begin{lemma}
\label{lemma:lyapunov-norm}

There exists a radius \(R > 0\)
such that, for all \(j \in [l+1]\),
the following properties hold.
\begin{enumerate}
    \item If \(x \in Y_j(\xi)\), \(\norm{A^{(N)}(\xi)x} \lessapprox e^{\chi_{j+1} N} \norm{x}\)
    as \(N \to \infty\), uniformly for \(x\) and \(\xi \in B_R\).
    \item If \(x \in Y_j(\xi)^\perp\), \(\norm{A^{(N)}(\xi)x} \gtrapprox e^{\chi_j N} \norm{x}\)
    as \(N \to \infty\), uniformly for \(x\) and \(\xi \in B_R\).
    \item If \(j \in [l]\) and \(E_{j+1}\) has a basis of eigenvectors of \(A(0)\),
    then \textit{(i)} holds with \(\lesssim\) instead of \(\lessapprox\).
    \item If \(E_j\) has a basis of eigenvectors of \(A(0)\),
    then \textit{(ii)} holds with \(\gtrsim\) instead of \(\gtrapprox\).
\end{enumerate}
\end{lemma}
\begin{proof}
We write \(q := m_{(j)}\). Note that,
using the definition of the cocycle and Lemma \ref{lemma:cocycle-norm-lq},
we have \(\norm{A^{(N)\wedge q}(\xi)} \lessapprox e^{(\chi_1 + \cdots + \chi_j)N}\)
as \(N \to \infty\) uniformly for \(\xi \in B_R\): we will use this fact multiple times in the proof.
\begin{enumerate}
    \item Let \(x \in Y_j(\xi)\).
    Then, by Lemma \ref{lemma:orthogonality-before},
    we have \[{\norm{A^{(N)}(\xi)x} \lesssim \norm{(1 - P_{(j)})A^{(N)}(\xi)x}}\]
    as \(N \to \infty\) uniformly for all small enough \(\xi\).
    Furthermore
    we have \[
        \abs{\langle A^{(N)\wedge q}(\xi)\beta_q(\xi), \alpha_q\rangle} \asymp \abs{\mu_{(q)}}^N
        \asymp \norm{A^{(N)\wedge q}(\xi)\beta_q(\xi)}
    \]
    as \(N \to \infty, \xi \to 0\).
    Therefore, if we let \(\pi_q^{(N)}(\xi)\)
    be the orthogonal projection onto \(\linspan A^{(N)\wedge q}(\xi)\beta_q(\xi)\)
    and use \Cref{lemma:projection-blade-bound},
    we can conclude that
    there exist some radius \(R > 0\),
    some number \(N_0 \in \bN\)
    and some constant \(0 < C < 1\)
    such that \(\norm{P_{(j)} - \pi_q^{(N)}(\xi)} \leq C\)
    for all \(\xi \in B_R\) and \(N\) large enough.
    Therefore, there exists some radius \(R > 0\)
    such that
    \[
        \norm{A^{(N)}(\xi)x} \lesssim \norm{(1 - P_{(j)})A^{(N)}(\xi)x}
         \lesssim \norm{(1-\pi_q^{(N)}(\xi))A^{(N)}(\xi)x}
    \]
    as \(N \to \infty\) uniformly for \(x\) and \(\xi \in B_R\).
    The last term is, by definition, orthogonal to \(A^{(N)\wedge q}\beta_q(\xi)\),
    therefore its norm can be estimated using \Cref{lemma:compose-blade} and \Cref{lemma:cocycle-norm-lq} to obtain \[
        \norm{(1-\pi_q^{(N)}(\xi))A^{(N)}(\xi)x}
        = \frac{\norm{A^{(N)\wedge (q+1)} (\beta_q(\xi) \wedge x)}}{{\norm{A^{(N)\wedge q} \beta_q(\xi)}}}
        \lessapprox e^{\chi_{j+1}N} \norm{x}
    \]
    as \(N \to \infty\) uniformly for \(x\) and \(\xi \in B_R\).
    \item Let \(x \in Y_j(\xi)^\perp\):
    that is, \(x \in \linspan \beta_q(\xi)\). 
    This means there exists a unique \((q-1)\)-blade \(\gamma(\xi,x)\) orthogonal to \(x\)
    such that \(\beta_q(\xi) = \gamma(\xi,x) \wedge x\).
    Therefore, there exists a radius \(R > 0\) such that \[
        \norm{A^{(N)}(\xi)x} \geq \frac{\norm{A^{(N)\wedge q}(\xi)\beta_q(\xi)}}{\norm{A^{(N)\wedge(q-1)}(\xi)\gamma(\xi,x)}}
        \gtrapprox e^{\chi_jN} \norm{x}
    \] as \(N \to \infty\) uniformly for \(x\) and \(\xi \in B_R\). Again,
    we are using \Cref{lemma:cocycle-norm-lq} in the last step.
    \item Assume \(E_{j+1}\) has a basis of eigenvectors of \(A(0)\).
    Then, applying \Cref{lemma:exterior-power-generalized-eigenvalues}
    to the basis \(v_1, \cdots, v_\ell\),
    we deduce that the top eigenspaces of \(A(0)^{\wedge (q+1)}\)
    also have a basis of eigenvectors,
    therefore \(\norm{A(0)^{\wedge (q+1)}} = \rho(A(0)^{\wedge (q+1)})\).
    Then, by \Cref{lemma:cocycle-norm-lq},
    \(\norm{A^{(N)\wedge(q+1)}(\xi)} \lesssim \rho(A(0)^{\wedge (q+1)})^N\)
    as \(N \to \infty\) uniformly for \(\xi \in B_R\).
    Then the same argument as in \textit{(i)} holds using \(\lesssim\) instead of \(\lessapprox\).
    \item Same as in \textit{(iii)}. \qedhere
\end{enumerate}
\end{proof}

\begin{lemma}\label{lemma:orthogonality-after}
Let \(j, k \in [l+1]\) with \(j < k\).

Then \[
    \norm{P_k A^{(N)}(\xi)P_j(\xi)x} \lessapprox \max(e^{\chi_k - \chi_j}, \lambda^{-1})^N \norm{A^{(N)}(\xi)P_j(\xi)x}
\]
as \(N \to \infty\)
uniformly for \(x \in \bC^\ell\) and \(\xi \in B_R\).
\end{lemma}
\begin{proof}[Sketch]
The following argument is inspired by \citet*{Dubiner1991}.
The intuition behind it is as follows:
if \(x \in Y_{j+1}(0)\),
multiplying it by \(A(0)\)
amplifies its component in \(E_j(0)\)
relative to the norm of \(x\),
so \(\frac{A(0)^Nx}{\norm{A(0)^Nx}}\) will \enquote{converge to \(E_j\)}.

Now, if \(x \in Y_{j+1}(\xi)\),
multiplying it by \(A(\xi)\)
will also tend to amplify the component in \(E_j(\xi)\),
but there is also an error of order \(O(\xi \norm{x})\) relative to the \(\xi = 0\) case:
in particular,
\(A(\xi)x\) may fail to be in \(E_j(\lambda^{-1}\xi)\) even if \(x \in E_j(\xi)\).
We will need to control this error, which accounts for the worse bound in the lemma,
when compared to \Cref{lemma:orthogonality-before}.
\end{proof}
\begin{proof}
Let \(y = P_j(\xi)x\): without loss of generality
we can assume \(y \neq 0\).
Also let \begin{align*}
    C_N &:= \frac{\norm{P_k A^{(N)}(\xi)y}}{\norm{A^{(N)}(\xi)y}} \\
    q_N &:= \frac{\norm{A^{(N-1)}(\xi)y}}{\norm{A^{(N)}(\xi)y}}e^{\chi_j} \\
    \mu &:= e^{\chi_k - \chi_j}
\end{align*}

Also let \(L_1 > 0\) be such that
\(e^{-\chi_k} \norm{A(\lambda^{-N}\xi) - A(0)} \leq L_1 \lambda^{-N}\)
for all \(\xi \in B_R\).

We claim that the following inequality holds:
\begin{equation} \label{eq:induction-dubiner-argument}
    C_N \leq \left(\prod_{k=1}^{N} q_k\right)\mu^N C_0
        + L_1 \sum_{k=1}^{N} \left(\prod_{n=N-k+1}^N q_n\right) \lambda^{-N+k} \mu^{k-1}
\end{equation}

Assuming that this inequality holds, we can prove the lemma.
To see this,
note that \Cref{lemma:lyapunov-norm}
implies that \(\prod_{n=1}^N q_n \approx 1\)
as \(N \to \infty\) uniformly for \(\xi \in B_R, x \in \bC^\ell\):
therefore, there exists an increasing function \(\epsilon : \bN \to \bR_{\geq 0}\)
such that \(\epsilon(N)^{-1} \leq \prod_{n=1}^N q_n \leq \epsilon(N)\)
for all \(N \in \bN\), \(x \in \bC^\ell\) and \(\xi \in B_R\),
and \(\limsup_{N \to \infty} \frac1N \log \epsilon(N) \leq 0\).
Then we have
\begin{align*}
    C_N &\leq \left(\prod_{k=1}^{N} q_k\right)\mu^N C_0
        + L_1 \sum_{k=1}^{N} \left(\prod_{n=N-k+1}^N q_n\right) \lambda^{-N+k} \mu^{k-1} \\
        &\leq \mu^N \epsilon(N) C_0 + L_1\sum_{k=1}^N \epsilon(N)\epsilon(N - k) \lambda^{-N+k} \mu^{k-1} \\
        &\leq \max(\mu, \lambda^{-1})^N (\epsilon(N)C_0 + L_1 N \epsilon(N)^2 \mu^{-1}) \\
        &\lessapprox \max(\mu, \lambda^{-1})^N
\end{align*}
as \(N \to \infty\) uniformly for \(x \in \bC^\ell, \xi \in B_R\),
which is the statement of the lemma.
Therefore, as long as we can prove \eqref{eq:induction-dubiner-argument},
we are done.

In order to prove \eqref{eq:induction-dubiner-argument},
note that \(C_N\) satisfies the following recursive inequality for \(N \in \bN\).
\begin{align*}
    C_{N} = \frac{\norm{P_k A^{(N)}(\xi)y}}{\norm{A^{(N)}(\xi)y}} &\leq q_N e^{-\chi_j}\frac{\norm{P_k A^{(N)}(\xi)y}}{\norm{A^{(N-1)}(\xi)y}} \\
&\leq q_N e^{-\chi_j}\frac{\norm{P_k A(0) A^{(N-1)}(\xi)y}}{\norm{A^{(N-1)}(\xi)y}}\\
            &\quad+ q_N e^{-\chi_j}\frac{\norm{P_k (A(\lambda^{-N+1}\xi) - A(0)) A^{(N-1)}(\xi)y}}{\norm{A^{(N-1)}(\xi)y}} \\
        &\leq q_{N} \mu C_{N-1} + L_1 q_{N}\lambda^{-N+1}
\end{align*}

Using this,
we can prove \eqref{eq:induction-dubiner-argument} by induction.
For \(N = 0\), the inequality reads \(C_0 \leq C_0\),
so nothing is to be done.
For higher \(N\) we have:
\begin{align*}
    C_N &\leq q_{N} \mu C_{N-1} + L_1 q_{N}\lambda^{-N+1} \\
        &\leq q_N\mu
            \left(\left(\prod_{n=1}^{N-1} q_n\right) \mu^{N-1} C_0 + \sum_{k=1}^{N-1} L_1 \left(\prod_{n=N-k}^{N-1} q_n\right) \lambda^{-N+k+1} \mu^{k-1}  \right)
            + L_1 q_{N}\lambda^{-N+1} \\
        &= \left(\left(\prod_{k=1}^{N} q_k\right)\mu^N C_0
        + L_1 \sum_{k=1}^{N-1} \left(\prod_{n=N-k}^N q_n\right) \lambda^{-N+k+1} \mu^{k}\right)
        + L_1 q_N \lambda^{-N+1} \\
        &= \left(\prod_{k=1}^{N} q_k\right)\mu^N C_0
        + L_1 \sum_{k=0}^{N-1} \left(\prod_{n=N-k}^N q_n\right) \lambda^{-N+k+1} \mu^{k} \\
        &= \left(\prod_{k=1}^{N} q_k\right)\mu^N C_0
        + L_1 \sum_{k=1}^{N} \left(\prod_{n=N-k+1}^N q_n\right) \lambda^{-N+k} \mu^{k-1}
\end{align*}
Which finishes the proof by induction of Inequality \eqref{eq:induction-dubiner-argument},
and hence the lemma.
\end{proof}

Now we are ready to prove \Cref{thm:main-theorem-lyapunov}.

\begin{proof}[Proof of \Cref{thm:main-theorem-lyapunov}]
Let \(A\) be a matrix function
which is locally Lipschitz at \(0\).
Note that it suffices to prove that there exists some norm on \(\bC^\ell\)
such that the inequalities in Theorem \ref{thm:main-theorem-lyapunov} hold:
as all norms on \(\bC^\ell\) are equivalent,
this means the inequalities will hold for any norm.

Pick a descending Jordan basis \(v_1, \cdots, v_\ell\) of \(A(0)\)
and pick an inner product of \(\bC^\ell\)
which makes this basis orthonormal.
Then we can define the projections \(P_j(\xi)\)
as in \Cref{def:lyapunov-projections},
using this inner product,
and all the lemmas we have proven so far hold.

Then, using the norm induced by this inner product,
we check the inequalities.
\begin{itemize}
    \item Inequalities \eqref{eq:lyapunov-norm},
    \eqref{eq:lyapunov-norm-diagonalizable} and \eqref{eq:lyapunov-norm-mult-1}
    follow from \Cref{lemma:lyapunov-norm},
    using the fact that \(E_j(\xi) = Y_{j-1}(\xi) \cap Y_{j}(\xi)^\perp\).
    \item Inequality \eqref{eq:lyapunov-orth-backward}
    follows from \Cref{lemma:orthogonality-after}.
    \item Inequality \eqref{eq:lyapunov-orth-forward}
    follows from \Cref{lemma:orthogonality-before}.
\end{itemize}
This concludes the proof.\end{proof}

\section{Main Argument}\label{chap:main-argument}

In this section,
we use the tools we have developed
to study the diffraction of substitution rules around the origin
and provide sufficient conditions for their hyperuniformity
in terms of their substitution matrix.

\subsection{Statement}
\label{sec:main-argument-statement}

Let \(\cS\) be a substitution rule on \(\bR^d\)
with \(\ell\) colors
and \(\bm{\Phi}\) a self-similar point process with rule \(\cS\).
Let \(\lambda > 1\) be the scaling constant
and \(\Mfull\) be the full substitution matrix
as defined in \Cref{chap:substitutions}.
We will assume that \(\cS\) is \emph{primitive}, i.e. \(\Mfull\) is primitive:
we also assume its Perron--Frobenius eigenvalue is \(\lambda^d\).

For each \(\mu \in \Spec \Mfull\),
denote the corresponding generalized eigenspace
by \(E_\mu\),
and define the following constants:
\begin{align*}
    \beta_\parallel(w) &:= d - \max\left\{ \log_\lambda \abs{\mu} \;\middle|\; \mu \in \Spec \Mfull \setminus \{\lambda^d\},\; w \not\in E_\mu^\perp \right\}, \\
    \beta_\perp(w) &:= d+1 - \max\left\{ \log_\lambda \abs{\mu} \;\middle|\; \mu \in \Spec \Mfull \setminus \{\lambda^d\},\; w \in E_\mu^\perp \right\}, \\
    \beta(w) &:= \min\left(\beta_\perp(w), \beta_\parallel(w)\right).
\end{align*}

Note these expressions make sense
even if \(0 \in \Spec \Mfull\),
in which case \(\log 0 = -\infty\):
we never need to subtract two infinities.
We also define \(\max \emptyset = -\infty\).

For every \(w \in \bC^\ell\),
let \(\hat{\eta}_w\) be the diffraction measure of \(\bm{\Phi}\)
with weights \(w\).
The most general theorem on the decay of \(\hat{\eta}_w\) around the origin we prove
is the following:

Recall that we defined the asymptotic notation \(\lessapprox\),
such that \(f(r) \lessapprox r^\alpha\) as \(r \to 0\)
if and only if for every \(\epsilon > 0\)
we have \(f(r) \lesssim r^{\alpha-\epsilon}\)
(see \Cref{sec:lyapunov-asymptotic-notation} for details).
\begin{theorem}
\label{thm:main-theorem}
The diffraction measure \(\hat{\eta}_w\) satisfies \[
    \hat{\eta}_w(B_r) \lessapprox r^{2 \beta(w)}
\] as \(r \to 0\).
\end{theorem}
For \(w = (1, \ldots, 1)^\top\),
we obtain \Cref{intro-origin-bound} from the introduction.

From \Cref{thm:main-theorem}, we get the following criterion for hyperuniformity.
\begin{corollary}
\label{cor:hyperuniformity-criterion}
Assume \(d \in \{1,2\}\).
Then, for every \(w \in \bC^\ell\),
the following holds:
if \(\abs{\mu} < \lambda^{\frac{d}{2}}\)
for all \(\mu \in \Spec \Mfull \setminus \{\lambda^d\}\) such that \(w \not\in E_\mu^\perp\),
then \(\bm{\Phi}\) is hyperuniform for weights \(w\).
\end{corollary}
\begin{proof}
By definition,
\(\bm{\Phi}\) is hyperuniform for weights \(w\)
if and only if \(\hat{\eta}_w(B_r) = o(r^d)\)
as \(r \to 0\),
where \(\hat{\eta}_w\) is the diffraction measure of \(\langle\bm{\Phi},w\rangle\).
Then we have:
\begin{itemize}
    \item For all \(\mu \in \Spec \Mfull \setminus \{\lambda^d\}\)
        with \(w \not\in E_\mu^\perp\),
        we have \(\abs{\mu} < \lambda^{\frac{d}2}\) by assumption.
        Therefore \(\beta_\parallel(w) > d - \frac{d}2 = \frac{d}{2}\).
    \item For all \(\mu \in \Spec \Mfull \setminus \{\lambda^d\}\)
    with \(w \in E_\mu^\perp\),
    we have \(\abs{\mu} < \lambda^d\), as \(\lambda^d\) is the PF eigenvalue of the substitution matrix.
        Therefore \(\beta_\perp(w) > d + 1 - d = 1\).
\end{itemize}
By Theorem \ref{thm:main-theorem},
we have \(\hat{\eta}_w(B_r) \lessapprox r^{2\beta(w)}\)
as \(r \to 0\).
But for \(d = 1, 2\), we have \(\beta(w) = \min (\beta_\perp(w), \beta_\parallel(w)) > \min(\frac{d}{2},1) \geq \frac{d}{2}\),
so \(\hat{\eta}_w(B_r) = o(r^d)\)
as \(r \to 0\).
This means \(\bm{\Phi}\) is hyperuniform for weights \(w\).
\end{proof}

\subsection{Proof of \Cref{thm:main-theorem}}
\label{sec:main-argument-proof}

Let \(\bm{A}\) be the normalized Fourier matrix of the substitution rule \(\cS\).
Let \(\chi_1 > \cdots > \chi_l > \chi_{l+1}\)
be the Lyapunov exponents of \(\bm{A}(0)\)
and \(E_1, \ldots, E_{l+1}\) their corresponding spaces
of generalized eigenvectors
with projections \(P_1, \ldots, P_{l+1}\),
as defined in \Cref{sec:lyapunov-exponents-converging-sequences}.
Also let \(R > 0\) be small enough so that the conclusions of Theorem \ref{thm:main-theorem-lyapunov}
hold for the normalized Fourier matrix \(\bm{A}\):
let \(P_1(\xi), \ldots, P_{l+1}(\xi)\)
be the projections from the theorem.
By \Cref{thm:renormalisation-measure-existence},
the diffraction matrix
\(\hat{\bm{H}}\)
of \(\bm{\Phi}\)
admits a self-similar density
\((\bm{h}, \nu)\) on \(B_R^\times\).

Fix \(w \in \bC^\ell\),
and define the two following subsets of \([l+1]\):
\begin{align*}
    \cJ_{\parallel}(w) &= \{j \in [l+1] \mid w \not\in E_j^\perp\} \\
    \cJ_{\perp}(w) &= \{j \in [l+1] \mid w \in E_j^\perp\}
\end{align*}

Then the constants from \Cref{thm:main-theorem}
can be restated as follows:
\begin{lemma}
\label{beta-with-lyapunov}

The constants \(\beta_\parallel\) and \(\beta_\perp\) from Theorem \ref{thm:main-theorem}
are given by
\begin{itemize}
    \item \(\beta_\parallel(w) = \min\{-\frac{\chi_j}{\log \lambda} \mid  j \in \cJ_{\parallel}(w) \setminus \{1\}\}\)
    \item \(\beta_\perp(w) = 1+\min\{-\frac{\chi_j}{\log \lambda} \mid  j \in \cJ_{\perp}(w)\setminus \{1\}\}\)
\end{itemize}
\end{lemma}
\begin{proof}
Recall that \(\beta_\parallel(w)\) is defined as
\[
    d - \max\left\{ \log_\lambda \abs{\mu} \;\middle|\; \mu \in \Spec \Mfull \setminus \{\lambda^d\},\; w \not\in E_\mu^\perp \right\}
\]

Using the fact that \(\bm{A}(0) = \lambda^{-d} \Mfull\),
we see that  \(\mu \in \Spec \Mfull\) satisfies \(\log \frac{\abs{\mu}}{\lambda^d} = \chi_j\)
if and only if generalized eigenspace \(E_\mu\) is contained in \(E_j\).
Therefore \(w \in E_j^\perp\)
holds if and only if \(w \in E_\mu^\perp\)
for every \(\mu \in \Spec \Mfull\) such that \(\log \frac{\abs{\mu}}{\lambda^d} = \chi_j\).
If this is the case, we have:
\[-\frac{\chi_j}{\log \lambda} = - \frac{\log \abs{\mu} - d \log \lambda}{\log \lambda} = d - \log_\lambda \abs{\mu}.\]

Therefore the definitions of \(\beta_\parallel(w)\) and \(\beta_\perp(w)\)
coincide with the equations in the statement of the lemma.
This concludes the proof.
\end{proof}

The proof of Theorem \ref{thm:main-theorem}
hinges on the following lemma.
\begin{lemma}\label{lemma:orthogonality-estimate-invariant}
For all \(j \in [l+1] \setminus \{1\}\)
and \(k \in \cJ_\parallel(w)\),
we have \[
    \norm{P_k \bm{A}^{(N)}(\xi)P_j(\xi)x} \lessapprox \lambda^{-\beta(w) N}\norm{P_j(\xi)x}
\]
as \(N \to \infty\) uniformly for \(\xi \in B_R\), \(x \in \bC^\ell\).
\end{lemma}
\begin{proof}
Depending on \(j\) and \(k\),
we will find \(\alpha, \gamma\) such that \(\alpha + \gamma \geq \beta(w)\) and the following inequalities hold:
\begin{equation} \label{eq:orthogonality-estimate}
    \begin{split}
    \norm{\bm{A}^{(N)}(\xi)P_j(\xi)x} \lessapprox \lambda^{-\alpha N}\norm{P_j(\xi)x} \\
    \norm{P_k 
    \bm{A}^{(N)}(\xi)P_j(\xi)x} \lessapprox \lambda^{-\gamma N}\norm{\bm{A}^{(N)}(\xi)P_j(\xi)x}
    \end{split}
\end{equation}
as \(N \to \infty\)
uniformly for \(\xi \in B_R\), \(x \in \bC^\ell\).
If we achieve this, we have \[
   \norm{P_k \bm{A}^{(N)}(\xi)P_j(\xi)x}
    =  \frac{\norm{P_k \bm{A}^{(N)}(\xi)P_j(\xi)x}}{\norm{\bm{A}^{(N)}(\xi)P_j(\xi)x}}\norm{\bm{A}^{(N)}(\xi)P_j(\xi)x}
    \lessapprox \lambda^{-(\alpha + \gamma)N}\norm{P_j(\xi)x}
\]
as \(N \to \infty\) uniformly for \(\xi \in B_R\), \(x \in \bC^\ell\),
 which concludes the proof of the lemma.
 (If the denominator of the fraction is \(0\) the left hand side is \(0\) as well,
 so the inequality still holds.)

In order to find \(\alpha, \gamma\),
we need to consider different cases.
We use \Cref{thm:main-theorem-lyapunov} for all of these,
using the characterization of \(\beta_\parallel\) and \(\beta_\perp\)
from \Cref{beta-with-lyapunov}.
\begin{itemize}
    \item \textbf{Case \(j \in \cJ_\parallel(w)\):} 
        In this case, Equation \eqref{eq:lyapunov-norm}
        from \Cref{thm:main-theorem-lyapunov}
        reads \[
            \norm{\bm{A}^{(N)}(\xi)P_j(\xi)x} \lessapprox e^{\chi_j N} \norm{P_j(\xi)x} 
        \]
        Furthermore,
        we have \[
        \norm{P_k \bm{A}^{(N)}(\xi)P_j(\xi)x} \lessapprox \norm{\bm{A}^{(N)}(\xi)P_j(\xi)x}
    \]
    Then, using the fact that \(e^{\chi_j N} = \lambda^{\frac{\chi_j}{\log \lambda} N}\),
        we see that
        Inequalities \eqref{eq:orthogonality-estimate}
        hold with \(\alpha = -\frac{\chi_j}{\log\lambda}\)
        and \(\gamma = 0\).
    As \(\alpha \geq \beta_\parallel(w) \geq \beta(w)\),
        we have \(\alpha + \gamma \geq \beta(w)\).
    \item \textbf{Case \(j \in \cJ_\perp(w), j > k\):}
        In this case,
        Equations \eqref{eq:lyapunov-norm} and \eqref{eq:lyapunov-orth-forward}
        from \Cref{thm:main-theorem-lyapunov} read:
        \begin{align*}
            \norm{\bm{A}^{(N)}(\xi)P_j(\xi)x} &\lessapprox e^{\chi_j N} \norm{P_j(\xi)x} \\
        \norm{P_k \bm{A}^{(N)}(\xi)P_j(\xi)x} &\lesssim \lambda^{-N} \norm{\bm{A}^{(N)}(\xi)P_j(\xi)x}
        \end{align*}
        Therefore, Inequalities \eqref{eq:orthogonality-estimate} hold
        if we set \(\alpha = -\frac{\chi_j}{\log\lambda}\) and \(\gamma = 1\):
        In this case \(\alpha + \gamma = 1 - \frac{\chi_j}{\log\lambda} \geq \beta_\perp(w) \geq \beta(w)\).
    \item \textbf{Case \(j \in \cJ_\perp(w), j < k\):}
        In this case,
            Equations \eqref{eq:lyapunov-norm} and \eqref{eq:lyapunov-orth-backward}
        from \Cref{thm:main-theorem-lyapunov} read:
        \begin{align*}
            \norm{\bm{A}^{(N)}(\xi)P_j(\xi)x} &\lessapprox e^{\chi_j N} \norm{P_j(\xi)x} \\
        \norm{P_k \bm{A}^{(N)}(\xi)P_j(\xi)x} &\lessapprox \max(e^{\chi_k - \chi_j}, \lambda^{-1})^N \norm{\bm{A}^{(N)}(\xi)P_j(\xi)x}
        \end{align*}
        Therefore, Inequalities \eqref{eq:orthogonality-estimate} hold
        if we set \(\alpha = -\frac{\chi_j}{\log\lambda}\) and \(\gamma = \min(\frac{\chi_j - \chi_k}{\log\lambda}, 1)\).
        In this case we have either \(\alpha + \gamma = -\frac{\chi_j}{\log\lambda} + 1 \geq \beta_\perp(w)\)
        or \(\alpha + \gamma = -\frac{\chi_j}{\log\lambda}+\frac{\chi_j - \chi_k}{\log\lambda} = -\frac{\chi_k}{\log\lambda} \geq \beta_\parallel(w)\).
        In both cases we have \(\alpha + \gamma \geq \beta(w)\).
\end{itemize}
Then, for every \(j, k\) as in the statement of the lemma,
we have found \(\alpha, \gamma\)
such that Inequalities \eqref{eq:orthogonality-estimate} hold
and \(\alpha + \gamma \geq \beta(w)\). This concludes the proof.
\end{proof}

In \Cref{lemma:orthogonality-estimate-invariant},
we exclude the case \(j = 1\).
This is justified by the following result.

\begin{lemma}\label{lemma:vanishing-first-space}
For \(\nu\)-almost every \(\xi \in B_R\),
we have \(P_1(\lambda^{-1}\xi)\bm{h}(\xi)P_1(\lambda^{-1}\xi)^* = 0\).
\end{lemma}
\begin{proof}
The first Lyapunov exponent of \(\bm{A}(0)\),
\(\chi_1 = 0\),
has multiplicity \(m_1 = 1\).
Therefore,
by Inequalities \eqref{eq:lyapunov-norm-mult-1} and \eqref{eq:lyapunov-orth-backward} 
of Theorem \ref{thm:main-theorem-lyapunov},
we have
\begin{align*}
    \norm{P_1(\xi)x} &\lesssim \norm{\bm{A}^{(N)}(\xi)P_1(\xi)x} \\
    &\leq \norm{P_1(\lambda^{-N}\xi)\bm{A}^{(N)}(\xi)P_1(\xi)x} + \sum_{j=2}^{l+1} \norm{P_j(\lambda^{-N}\xi)\bm{A}^{(N)}(\xi)P_1(\xi)x} \\
    &\lesssim \norm{P_1(\lambda^{-N}\xi)\bm{A}^{(N)}(\xi)P_1(\xi)x}
\end{align*}
as \(N \to \infty\)
uniformly for all \(\xi \in B_R, x \in \bC^\ell\). 
This means there exist some \(c > 0\) and \(N_0 \in \bN\)
such that \[
    \norm{P_1(\lambda^{-N}\xi)\bm{A}^{(N)}(\xi)P_1(\xi)x} \geq c\norm{P_1(\xi)x}
\]
for all \(\xi \in B_{R}\), \(N \geq N_0\).

Furthermore, recall that \(A^{(N)}(\xi)Y_1(\xi) \subset Y_1(\lambda^{-N}\xi)\).
By definition, \(Y_1(\xi)\) is the image of \(I - P_1(\xi)\),
so \(\bm{A}^{(N)}(\xi)(I - P_1(\xi)) = (I - P_1(\lambda^{-N}\xi))\bm{A}^{(N)}(\xi)(I - P_1(\xi))\)
and therefore \[P_1(\lambda^{-N}\xi)\bm{A}^{(N)}(\xi)P_1(\xi) = P_1(\lambda^{-N}\xi)\bm{A}^{(N)}(\xi).\]

As \(\bm{h}(\xi)\) is positive semidefinite for all \(\xi\),
we can find \(f_1(\xi), \ldots, f_\ell(\xi)\)
such that \(\bm{h}(\xi) = \sum_{t=1}^\ell f_t(\xi)f_t(\xi)^*\).
Therefore \begin{align*}
    \norm{\bm{h}(\lambda^{-N}\xi)} &= \norm{\bm{A}^{(N)}(\lambda^{-1}\xi)\bm{h}(\xi)\bm{A}^{(N)*}(\lambda^{-1}\xi)} \\
    &\geq \norm{P_1(\lambda^{-N-1}\xi)\bm{A}^{(N)}(\lambda^{-1}\xi)\bm{h}(\xi)\bm{A}^{(N)*}(\lambda^{-1}\xi)P_1(\lambda^{-N-1}\xi)^*} \\
    &= \sum_{t=1}^\ell \norm{P_1(\lambda^{-N-1}\xi)\bm{A}^{(N)}(\lambda^{-1}\xi)f_t(\xi)}^2 \tag{$*$} \\
    &= \sum_{t=1}^\ell \norm{P_1(\lambda^{-N-1}\xi)\bm{A}^{(N)}(\lambda^{-1}\xi)P_1(\lambda^{-1}\xi)f_t(\xi)}^2 \\
    &\geq c\sum_{t=1}^\ell \norm{P_1(\lambda^{-1}\xi)f_t(\xi)}^2 \\
    &= c\norm{P_1(\lambda^{-1}\xi)\bm{h}(\xi)P_1(\lambda^{-1}\xi)^*}
\end{align*}
for all \(\xi \in B_{R}\), \(N \geq N_0\).
Note that the equality \((*)\) holds because \(P_1(\lambda^{-N-1}\xi)\) is a \emph{rank \(1\)} projection.

On the other hand,
the diffraction matrix satisfies \(\hat{\bm{H}}(\{0\}) = 0\) \citep*[Lemma 2.6]{hartnick2022hyperuniformity}
so \(\bm{h}(\lambda^{-N}\xi) \to 0\) as \(N \to \infty\).
Then the inequality we just proved implies \(\norm{P_1(\lambda^{-1}\xi)\bm{h}(\xi)P_1(\lambda^{-1}\xi)^*} = 0\),
as we wanted to show.
\end{proof}

\begin{proof}[Proof of Theorem \ref{thm:main-theorem}]
By Lemma \ref{lemma:uniformity-from-density}, it suffices to show
\begin{equation}\label{eq:goal}
    \langle \bm{h}(\lambda^{-N}\xi)\,w, w\rangle
    \;\lessapprox\;
    \lambda^{-2\beta(w)N}\,\norm{\bm{h}(\xi)}
    \qquad\text{as }N\to\infty,
\end{equation}
uniformly for \(\xi\in B_R\).

As \(\hat{\bm{H}}(f)\) is a positive definite matrix for all positive \(f \in C_c(\bR^d)\),
one can show that the density \(\bm{h}(\xi)\)
is a positive semidefinite matrix for \(\nu\)-almost every \(\xi\).
Hence for \(\nu\)-almost every \(\xi\) there exist orthogonal vectors
\(f_1(\xi), \ldots, f_\ell(\xi)\) (possibly some zero) with
\[
    \bm{h}(\xi)=\sum_{t=1}^\ell f_t(\xi)f_t(\xi)^*.
\]
Therefore
\begin{equation*}\label{eq:decomp}
    \langle \bm{h}(\lambda^{-N}\xi)\,w, w\rangle
    = \sum_{t=1}^\ell \bigl|\langle \bm{A}^{(N)}(\lambda^{-1}\xi)f_t(\xi),\,w\rangle\bigr|^2,
\end{equation*}
and it is enough to prove
\begin{equation*}\label{eq:star}
    \bigl|\langle \bm{A}^{(N)}(\lambda^{-1}\xi)f_t(\xi),\,w\rangle\bigr|
    \;\lessapprox\;
    \lambda^{-\beta(w)N}\,\norm{f_t(\xi)},
\end{equation*}
as \(N\to\infty\), uniformly for \(\xi\in B_R\) and for all \(t\in[\ell]\).

We can decompose this further using the projections \(P_j(\lambda^{-1}\xi)\) and \(P_j\)
from Theorem \ref{thm:main-theorem-lyapunov}.
We have \(\sum_{j=1}^{l+1}P_j = \sum_{j=1}^{l+1}P_j(\lambda^{-1}\xi) = I\),
so we can bound
\begin{align*}
\left|\langle \bm{A}^{(N)}(\lambda^{-1}\xi)f_t(\xi), w\rangle\right|
&= \left|\left\langle \Big(\sum_{k=1}^{l+1}P_k\Big)\bm{A}^{(N)}(\lambda^{-1}\xi)
    \Big(\sum_{j=1}^{l+1}P_j(\lambda^{-1}\xi)f_t(\xi)\Big), w\right\rangle\right|\\
&= \left|\sum_{j,k=1}^{l+1}\left\langle \bm{A}^{(N)}(\lambda^{-1}\xi)P_j(\lambda^{-1}\xi)f_t(\xi), P_k^* w\right\rangle\right|\\
&\le \sum_{j,k=1}^{l+1}\left|\left\langle \bm{A}^{(N)}(\lambda^{-1}\xi)P_j(\lambda^{-1}\xi)f_t(\xi), P_k^* w\right\rangle\right|.
\end{align*}

Therefore, in order to prove the theorem,
it suffices to prove
\[
    \abs{\langle \bm{A}^{(N)}(\lambda^{-1}\xi)P_j(\lambda^{-1}\xi)f_t(\xi), P_k^*w\rangle}
    \lessapprox \lambda^{-\beta(w) N}\norm{P_j(\lambda^{-1}\xi)f_t(\xi)} \norm{P_k^*w}
    \tag{\(\star\)}
\]
as \(N \to \infty\)
uniformly for all \(\xi \in B_R\),
\(t \in [\ell]\) and \(j,k \in [l+1]\).
If we are able to prove this inequality,
we are done.

It suffices to prove \((\star)\)
for all \(\xi, t, j, k\)
such that \(P_j(\lambda^{-1}\xi)f_t(\xi) \neq 0\)
and \(P_k^* w \neq 0\),
as otherwise both sides of the asymptotic inequality are \(0\).
But:
\begin{itemize}
    \item By \Cref{lemma:vanishing-first-space},
we have \(P_1(\lambda^{-1}\xi)f_t(\xi) = 0\).
    \item If \(k \not\in \cJ_\parallel(w)\),
    we have \(k \in \cJ_\perp(w)\).
    By definition, this means \(w \perp E_k\):
    as \(P_k\) is a projection with image \(E_k\),
    this means \(P_k^* w = 0\)
\end{itemize}
This means that we only need to prove the inequality \((\star)\)
for the \(j, k \in [l+1]\)
such that \(j \neq 1\)
and \(k \in \cJ_\parallel(w)\).
These are precisely the indices
we considered in Lemma \ref{lemma:orthogonality-estimate-invariant}:
therefore, Lemma \ref{lemma:orthogonality-estimate-invariant}
proves \((\star)\)
uniformly for all \(\xi \in B_R\),
\(t \in [\ell]\) and \(j,k \in [l+1]\).
This completes the proof.
\end{proof}
\section{Construction of self-similar point processes}
\label{sec:construction}

In this section, we turn our attention to \emph{stone substitution rules},
and construct self-similar point processes associated to them.
Let \(\cS\) be a primitive substitution rule on \(\bR^d\) with \(\ell\) colors.

\begin{definition}
The \emph{substitution map}
\(\varrho : \cM(\bR^d) \to \cM(\bR^d)\) associated to \(\cS\)
is defined as follows:
if \(\bm{\mu} = (\mu_j)_{j \in [\ell]} \in \cM(\bR^d)^\ell\),
\(\varrho(\bm{\mu})\) has the entries \[
    \varrho(\bm{\mu})_k := \sum_{j=1}^\ell \sum_{x \in \Delta_{kj}} T_x D_\lambda \mu_j, \quad \forall k \in [\ell].
\]
\end{definition}

\begin{definition}
Let \(\bm{\mu} = (\mu_j)_{j\in [\ell]}\)
be a \(\bC^\ell\)-measure
and \(\tau_1, \ldots, \tau_\ell, S \subset \bR^d\)
be closed subsets.
We say \(\bm{\mu}\) \emph{tiles \(S\) with prototiles \(\tau_1, \ldots, \tau_\ell\)}
if we have \[
    \bone_S = \sum_{j=1}^\ell \mu_j * \bone_{\tau_j}
\] as an equality of \(L^\infty\) functions.
In this case, the set \[
    \{T_x \tau_j \mid j \in [\ell], x \in \supp \mu_j\}
\]
is a \emph{tiling of \(S\) with prototiles \(\tau_1, \ldots, \tau_\ell\)}.
\end{definition}

For \(j \in [\ell]\), let \(\bm{o}_j \in \cM(\bR^d)^\ell\)
be the vector measure with \((\bm{o}_j)_k = \delta_0\) if \(k = j\):
we can think of \(\bm{o}_j\) as a \enquote{colored point}.
\begin{definition}
\index[defs]{stone}
Let \(\cS\) be a primitive substitution rule on \(\bR^d\).
We say \(\cS\) is \emph{stone}
if the following holds:
there exist compact sets \({\tau_1, \ldots, \tau_\ell \subset \bR^d}\)
with positive volume
such that, for all \(j \in [\ell]\),
\(\varrho(\bm{o}_j)\) tiles \(D_\lambda \tau_j\)
with prototiles \(\tau_1, \ldots, \tau_\ell\).
\end{definition}

We say the compact sets \(\tau_1, \ldots, \tau_\ell\)
are the \emph{canonical prototiles} of \(\cS\):
these are uniquely determined by \(\cS\)
\citep*{lagarias2002substitutiondelonesets}.
Note that, if \(\cS\) is stone, the PF eigenvalue of its substitution matrix
is necessarily equal to \(\lambda^d\),
as the volumes of the prototiles form a positive eigenvector.

If \(\cS\) is stone, the following theorem holds:
this result is a restatement of well-known facts about stone substitution rules \citep*{Robinson2004, Baake2013, FRANK2008295}.
\begin{theorem}
Let \(\cS\) be a primitive, stone substitution rule with canonical prototiles \(\tau_1, \ldots, \tau_\ell\).
Then there exists a unique ergodic self-similar point process \(\bm{\Phi} = (\Phi_j)_{j\in[\ell]}\).
We say \(\bm{\Phi}\) is the \emph{vector point process associated to \(\cS\)}.
\end{theorem}
\begin{proof}
As is well known in the literature of aperiodic order,
given a primitive stone substitution rule \(\cS\),
one can always find a \emph{self-similar tiling},
i.e. a vector point measure \(\bm{\Lambda}\)
such that \begin{itemize}
    \item \(\bm{\Lambda}\) tiles \(\bR^d\) with prototiles \(\tau_1, \ldots, \tau_\ell\), and
    \item there exists \(k \in \bN\) such that \(\varrho^k(\bm{\Lambda}) = \bm{\Lambda}\).
\end{itemize}
Then one can define the \emph{substitution space} \[
    \Omega_\cS := \overline{\{T_x \bm{\Lambda} \mid x \in \bR^d\}}
\]
using the Chabauty--Fell topology on closed subsets of \(\bR^d\),
also called the \enquote{local rubber topology} in the field of aperiodic order.
In the case where \(\bm{\Lambda}\) has finite local complexity (FLC),
this is equivalent to the local topology on Delone sets,
which is the case usually considered in the literature
\citep*{Robinson2004, Baake2013, FRANK2008295,solomyak-notes}.
The general case is due to \citet*{LeeNonFLC}:
see also \citet*{Robinson2004}.

The space \(\Omega_\cS\) is compact and uniquely ergodic.
The unique probability measure on \(\Omega_\cS\)
is also invariant under the substitution map \(\varrho\),
and hence defines a self-similar point process \(\bm{\Phi}\)
with rule \(\cS\).

For the uniqueness,
let \(\bm{\Phi}'\) be an ergodic self-similar point process with rule \(\cS\),
and let \(\bP'\) be its distribution.
Using its self-similarity, one can show that the support of \(\bP'\)
contains the orbit of a self-similar tiling \(\bm{\Lambda}\),
and its closure has nonzero measure.
Therefore, by the unique ergodicity, \(\bm{\Phi}'\) must be the same as \(\bm{\Phi}\).
\end{proof}

Furthermore, by construction, \(\bm{\Phi}\) defines tilings of \(\bR^d\),
as follows.
\begin{theorem}
Let \(\bm{\Phi} = (\Phi_j)_{j \in [\ell]}\) be the vector point process
associated to \(\cS\),
and let \(\tau_1, \ldots, \tau_\ell\) be the prototiles of \(\cS\).
Then \(\bm{\Phi}\) almost surely tiles \(\bR^d\) with prototiles \(\tau_1, \ldots, \tau_\ell\).
\end{theorem}

The tilings \(\{T_x \tau_j \mid j \in [\ell], x \in \Phi_j\}\)
are called \emph{substitution tilings with rule \(\cS\)}.

\begin{remark}
In the literature on aperiodic order,
one often considers tilings instead of point sets as the main objects of study.
These viewpoints are essentially equivalent:
given a substitution tiling, one can obtain a locally finite set
by marking a point in each canonical prototile
in a consistent way:
then the set of marked points of each color
may be interpreted as determining the relative locations of the tiles of that color.
The particular point set obtained from a tiling depends on this choice,
but the inflation matrix remains unchanged:
therefore, our bounds for the diffraction and our criteria for hyperuniformity
do not depend on this choice.
\end{remark}

\section{Examples}
\label{sec:examples}
Now we apply the results of \Cref{chap:main-argument}
to the self-similar point processes constructed in the previous section,
and discuss the results.

\subsection{One-dimensional substitution rules with two colors}
\label{sec:main-argument-one-dimensional-case}

First,
we turn our attention to the simplest nontrivial substitution rules possible.
Let \(\cS\) be a primitive, stone substitution rule on \(\bR\) with \(\ell = 2\),
and let \(\bm{\Phi}\) be the vector point process with rule \(\cS\).
In this case,
the substitution matrix only has one non-PF eigenvalue:
Let \(\lambda\) be the scaling constant of \(\cS\),
\(\mu_2\) be the smallest eigenvalue of the substitution matrix \(\Mfull\),
and \(v_2\) be the corresponding right eigenvector.
Then \Cref{thm:main-theorem} specializes to the following:

\begin{corollary}
\label{corollary:main-theorem-2-colors}
For all \(w \in \bC^2\),
let \(\hat{\eta}_w\) be the diffraction of \(\bm{\Phi}\) with weights \(w\).
Then we have \[
    \hat{\eta}_w(B_r) \lessapprox r^{2 \beta(w)}
\]
as \(r \to 0\),
where \[
    \beta(w) := \begin{cases}
        2-\log_\lambda\abs{\mu_2} & \text{if } w \perp v_2 \\
        1-\log_\lambda\abs{\mu_2} & \text{otherwise}
    \end{cases}
    \]
Furthermore, the inequality holds with \(\lesssim\)
unless \(\mu_2 = 0\),
in which case we have \(\hat{\eta}_w(B_r) \lessapprox 0\)
as \(r \to 0\).
\end{corollary}

(The \(\lesssim\) case
does not immediately follow from \Cref{thm:main-theorem}
and requires a more careful analysis of the proof: see \citep*{thesis})

In the generic case where \(w\) is not orthogonal to the second eigenvector,
this is the bound originally conjectured by \citet*{Ouz2019}
and has been checked for several examples
by the same authors and also by \citet*{Baake_2019_origin},
as we see in the examples below.

In particular, our theorem proves the hyperuniformity
of certain already known specific examples:
see \citep*{thesis} for the definitions.
\begin{example}
\leavevmode
\begin{itemize}
    \item Let \(\cS\) be the Fibonacci substitution rule,
    which is the geometric version of the symbolic substitution \(1 \mapsto 12, 2 \mapsto 1\).
    Its scaling constant is \(\lambda = \tau = \frac{1+\sqrt{5}}{2}\)
    and its second eigenvalue is \(\mu_2 = -\lambda^{-1}\),
    hence we have \(\hat{\eta}_w(B_r) \lesssim r^4\) as \(r \to 0\) for all weights \(w\),
    as already shown by \citet*{Baake_2019_origin}.
    If \(w = (\tau^{-1}, \tau^{-2})^\top\),
    we have the improved bound \(\hat{\eta}_w(B_r) \lesssim r^6\):
    this is a new result.
    \item Let \(\cS\) be the Thue--Morse substitution rule,
    which is the geometric version of the symbolic substitution \(1 \mapsto 12, 2 \mapsto 21\).
    Its scaling constant is \(\lambda = 2\)
    and its second eigenvalue is \(\mu_2 = 0\),
    hence we have \(\hat{\eta}_w(B_r) \lessapprox 0\) as \(r \to 0\) for all weights \(w\):
    in other words, \(\hat{\eta}_w(B_r)\) decays faster than any polynomial as \(r \to 0\).
    Note that this result was already obtained by \citet*{thue-morse} with special methods,
    but in our case it follows from a general theorem.
\end{itemize}
\end{example}

\subsubsection{Non-invariance of hyperuniformity under MLD}
\label{subsec:non-mld-invariance}
Using our methods, we prove that
hyperuniformity is not an MLD invariant,
even among FLC sets
(see \citet*[Chapter 5]{Baake2013} for definitions of these terms).
This contradicts prevailing conjectures in the field \citep*{invariance-conjecture}.

\begin{example}\label{ex:anti-tilings}
Define the substitution rule \(\cS = (\lambda, \bm{\Delta})\) on \(\bR\)
by setting \(\lambda = 5\) and \[
    \bm{\Delta} = \begin{pmatrix}
        \{0, 1, 2, 3\} & \{-2-\frac12\} \\
        \{3+\frac12\} & \{-1, 0, 1, 2\}
    \end{pmatrix}
\]
This substitution rule is primitive and stone
with canonical prototiles \(\tau_1 = [0,1]\), \(\tau_2 = [-\frac12,\frac12]\)
and substitution matrix \(
    \Mfull = \begin{pmatrix}
        4 & 1 \\
        1 & 4
    \end{pmatrix}
\)
\end{example}

Let \(\bm{\Phi} = (\Phi_1, \Phi_2)\) be the vector point process associated to this substitution rule
and \(\hat{\eta}_w\) be its diffraction with weights \(w \in \bC^2\). We consider two special cases:
\begin{enumerate}
    \item Consider the weights \(w = (1,1)^\top\),
    such that \(\langle\bm{\Phi},w\rangle = \Phi_1 \dot{\cup} \Phi_2\).
    Then we have \(\hat{\eta}_w(B_r) \lesssim r^{2(2-\log_5 3)}\).
    In particular, because \(2(2-\log_5 3) =\, {\sim}2.6\),
    the vector point process associated to \(\cS\) is Class I hyperuniform for constant weights.
    \item Consider the weights \(w = (1,0)^\top\),
    such that \(\langle \bm{\Phi}, w\rangle = \Phi_1\).
    Then \Cref{corollary:main-theorem-2-colors}
    only implies \(\hat{\eta}_w(B_r) \lesssim r^{2(1-\log_5 3)}\),
    which is not enough to prove hyperuniformity.
\end{enumerate}
The second case suggests that \(\Phi_1\) may be non-hyperuniform.
Now we will prove that this is in fact the case.

\begin{theorem}\label{thm:anti-tilings-non-uniformity}
Let \(w = (1,0)^\top\).
Then \(\hat{\eta}_w(B_r) \gtrsim r^{\alpha}\) as \(r \to 0\),
where \(\alpha = 2 - 2\log_5 3\).
In particular, \(\Phi_1\) is not hyperuniform.
\end{theorem}

We prove this by finding an explicit self-similar set \(\bm{\Lambda}\)
which is \emph{limit-periodic}, i.e. a disjoint union of periodic sets.
This allows us to compute its diffraction explicitly.
See \citep*{limit-periodic}, \citep*[Sections 9.4.4, 9.4.5]{Baake2013}
for more examples of limit-periodic substitution sets and their diffraction.
We define the sequence of finite subsets \((S_n)_{n \in \bN_0}\) of \(\bR\)
by \begin{align*}
    S_0 &:= \{0\}, \\
    S_{n+1} &:= D_5 S_n \oplus \{1,2,3\}, n \in \bN_0,
\end{align*}
where by \(\oplus\) we mean the Minkowski sum of subsets of \(\bR\).
Then we define the pair \(\bm{\Lambda} = (\Lambda_1, \Lambda_2)\)
of point measures by \begin{align*}
    \Lambda_1 &:= \dot{\bigcup}_{n=0}^\infty \dot{\bigcup}_{x \in S_n} T_x D_{5^n} \bZ, \\
    \Lambda_2 &:= \dot{\bigcup}_{n=0}^\infty \dot{\bigcup}_{x \in S_n} T_{\frac12-x} D_{5^n} \bZ
\end{align*}

Then we have the following facts.
See \citep*[Chapter 9]{Baake2013} for the definition of the diffraction of a point set/measure,
and \citep*{BaakeBirknerGrimm2015,BaakeBirknerMoody2009} for its correspondence with the diffraction of the corresponding point process.
\begin{lemma}\leavevmode
\begin{enumerate}
    \item The pair \(\bm{\Lambda} = (\Lambda_1, \Lambda_2)\)
    satisfies \(\varrho(\bm{\Lambda}) = \bm{\Lambda}\).
    \item The diffraction of \(\bm{\Phi}\) with weights \(w\)
    is equal to the diffraction of \(w_1 \delta_{\Lambda_1} + w_2 \delta_{\Lambda_2}\),
    where \(\delta_{\Lambda_1}\) and \(\delta_{\Lambda_2}\)
    are the Dirac combs of \(\Lambda_1\) and \(\Lambda_2\) respectively.
\end{enumerate}
\end{lemma}

\begin{proof}[Proof of \Cref{thm:anti-tilings-non-uniformity}]
For \(w = (1,0)^\top\),
\(\hat{\eta}_w\) is the diffraction of the set \(\Lambda_1\).
By direct computation using the Poisson summation formula,
we have \[
    \hat{\eta}_w(\{5^{-n}\}) = \left|\sum_{k=n}^\infty \sum_{x \in S_k} e^{2 \pi i x 5^{-n}} 5^{-k}\right|^2 \gtrsim (\frac35)^{2n}.
\]
Furthermore, by the definition of \(\alpha\), \((\frac35)^2 = 5^{-\alpha}\),
so, by \Cref{lemma:uniformity-from-density},
this is enough to prove \(\hat{\eta}_w(B_r) \gtrsim r^{\alpha}\),
as desired.
\end{proof}

However, the following also holds.
\begin{theorem}\label{thm:anti-tilings-mld}
The sets \(\Lambda_1\) and \(\Lambda_1 \cup \Lambda_2\) are MLD equivalent FLC sets.
\end{theorem}

This means we have two MLD equivalent FLC sets
such that one gives rise to a hyperuniform point process and the other does not.

\begin{proof}
Given \(\Lambda_1\), one can reconstruct \(\Lambda_2\)
as the unique Delone set such that \((\Lambda_1, \Lambda_2)\) tiles \(\bR\)
with prototiles \(\tau_1, \tau_2\). This is clearly done locally.

Given \(\Lambda_1 \cup \Lambda_2\), one can recognize which points belong to \(\Lambda_1\)
by looking at the gaps between points.
Specifically:
\begin{itemize}
    \item Whenever two points are at distance \(1\), this means they have the same color,
    so either both or neither belong to \(\Lambda_1\).
    \item Whenever two points are at distance \(\frac12\), only the right one belongs to \(\Lambda_1\).
    \item Whenever two points are at distance \(\frac32\), only the left one belongs to \(\Lambda_1\).
\end{itemize}
As there are never more than five points with the same color in a row, this means that,
in order to tell whether a point belongs to \(\Lambda_1\),
it suffices to look at the next five points to the right:
in particular, it suffices to look at a uniformly bounded neighbourhood of the point.
\end{proof}

\subsection{Two-dimensional examples}\label{sec:2d-examples}
Now we turn our attention to two-dimensional substitution rules.
Here our methods greatly improve on existing results,
proving hyperuniformity for the vector point processes
associated to a wide class of substitution rules
where this was not previously known.
Virtually all known examples of two-dimensional substitution rules
have rotational or reflectional symmetry, which plays an essential role.

Let \(\cS\) be a primitive, stone substitution rule on \(\bR^d\)
with \(\ell\) colors,
\(\lambda > 1\) its scaling constant,
\(\bm{\Delta} = (\Delta_{jk})_{j,k \in [\ell]}\) be its displacement matrix,
\(\Mfull\) be its full substitution matrix
and \(\tau_1, \ldots, \tau_\ell\) its canonical prototiles.

\subsection{Symmetric substitution systems}
\label{sec:symmetric-substitution-systems}
First,
we define symmetric substitution systems,
and show how one can efficiently construct them and check their basic properties.
Let \(\cS\) be a primitive, stone substitution rule on \(\bR^d\).

\begin{definition}\label{def:symmetry}
\index[defs]{symmetric!substitution rule}
Let \(G < O(d)\) be a finite subgroup
equipped with an action on \([\ell]\),
in addition to its usual action on \(\bR^d\).
We say \(\cS\) is \emph{\(G\)-symmetric} if
it satisfies the following conditions:
\begin{align*}
\Delta_{g\cdot j,\, g\cdot k} &= g\cdot \Delta_{j,k} &&\text{for all } j,k\in[\ell],\ g\in G,\\
g\cdot \tau_j &= \tau_{g\cdot j} &&\text{for all } j\in[\ell],\ g\in G.
\end{align*}
\end{definition}
This action of \(G\) has important consequences for the associated point processes.
For us, the crucial fact is that,
if \(\cS\) is a \(G\)-symmetric substitution rule,
then the full substitution matrix \(\Mfull\)
commutes with the action of \(G\) on \(\bC^\ell\):
this simplifies the eigenvalues and eigenvectors of \(\Mfull\).
We state this in terms of the \emph{spherical substitution matrix},
which is a smaller matrix counting tiles only up to the action of \(G\).

\begin{definition}
\index[defs]{substitution matrix!spherical}
\label{def:spherical-substitution-matrix}
Let \(\cJ_G\) be a set of representatives of the orbits of \(G\) on \([\ell]\),
and let \(i_1 < \ldots< i_{\ell_0}\) be the elements of \(\cJ_G\).
Then the \emph{spherical substitution matrix} \(\Msph\)
is the \(\ell_0 \times \ell_0\) matrix defined by
\[    (\Msph)_{jk} = \sum_{g \in G} \# \Delta_{i_j, g \cdot i_k}
\]
for \(j,k \in [\ell_0]\).
\end{definition}

In the literature, both the full substitution matrix
(defined in \Cref{chap:substitutions})
and the spherical substitution matrix
may be referred to as \enquote{the substitution matrix of \(\cS\)}:
we always specify which one we mean from now on.
The spherical substitution matrix also depends on a choice of representatives \(\cJ_G\),
but this only alters the matrix by a permutation similarity,
and hence does not affect its spectral properties:
in particular, \(\Spec \Msph\) is independent of the choice of \(\cJ_G\).

Using \Cref{cor:hyperuniformity-criterion}
we can obtain a sufficient condition for hyperuniformity
in terms of the spherical substitution matrix alone.
In the examples, we see that determining the spherical substitution matrix
is much faster than determining the full substitution matrix.
\begin{corollary}
\label{corollary:symmetric-substitution-systems-uniformity}
Let \(\cS\) be a \(G\)-symmetric substitution rule
on \(\bR^2\).
Then the following holds:
if \(\abs{\mu} < \lambda\)
for all \(\mu \in \Spec \Msph \setminus \{\lambda^2\}\)
and \(w \in \bC^\ell\) is a \(G\)-invariant vector,
then the point process associated to \(\cS\) is hyperuniform for weights \(w\).
\end{corollary}
\begin{proof}
By the definition of a symmetric substitution rule,
one can check that the full substitution matrix \(\Mfull\)
commutes with the action of \(G\) on \(\bC^\ell\).
Hence, every generalized eigenspace of \(\Mfull\)
can be decomposed into a \(G\)-invariant part and an orthogonal part.
One can check that \(\Msph\) is the transformation matrix
of the restriction of \(\Mfull\) to the space of \(G\)-invariant vectors,
with a correct choice of basis.

Let \(\mu \in \Spec \Mfull\) be an eigenvalue such that \(w \not\perp E_\mu\).
By the above discussion, \(\mu\) must have a \(G\)-invariant eigenvector,
hence \(\mu \in \Spec \Msph\).
Therefore, when checking the hypotheses of \Cref{cor:hyperuniformity-criterion},
it suffices to consider eigenvalues in the spectrum of \(\Msph\).
\end{proof}

For every example we consider,
the symmetry group \(G\)
is either the cyclic group of rotations or the dihedral group.
\begin{itemize}
    \item \(G = C_n\), where the generator \(R\) acts on \(\bR^2\) by rotation by \(\frac{2\pi}{n}\) radians counterclockwise.
    For \(i \in [n]\), we set \(g_i := R^{i-1}\), so we have \(G = \{g_1, g_2, \ldots, g_n\}\).
    \item \(G = D_{n}\), the dihedral group of order \(2n\),
    generated by the rotation \(R\) by \(\frac{2\pi}{n}\) radians counterclockwise,
    and the reflection \(S\) in the \(x\)-axis.
    For \(i \in [n]\), we set \(g_i := R^{i-1}\) and \(g_{n+i} := R^{i-1}S\):
    this gives us an ordering \(G = \{g_1, g_2, \ldots, g_{2n}\}\).

\end{itemize}

\begin{remark}[A picture is worth a thousand words]
For these examples,
explicitly writing down the displacement matrix \(\bm{\Delta}\)
would often be cumbersome and not very enlightening.
As is common in the literature,
we define the substitution rule
by showing a picture of the patches defined by \(\varrho(\bm{o}_j)\)
for all \(j \in \cJ_G\),
where \(\cJ_G\) is a set of representatives of the orbits of \(G\) on \([\ell]\):
as long as one marks each tile in a way that makes its color clear,
this is enough to determine \(\varrho(\bm{o}_j)\)
for all \(j \in \cJ_G\),
and hence the substitution rule \(\cS\).
The way \(G\) operates on the colors
will also be apparent from the pictures.
\end{remark}

\subsection{Two-dimensional examples}
\label{sec:symmetric-substitution-systems-examples}
Now we apply \Cref{corollary:symmetric-substitution-systems-uniformity}
to known examples of substitution rules.
We give the definition of the chair rule in detail
and refer to the literature for the rest.
All of the results of this section are new.

\subsubsection{The chair substitution rule}

\begin{figure}[h]
        \centering
        \adjincludegraphics[width=\linewidth]{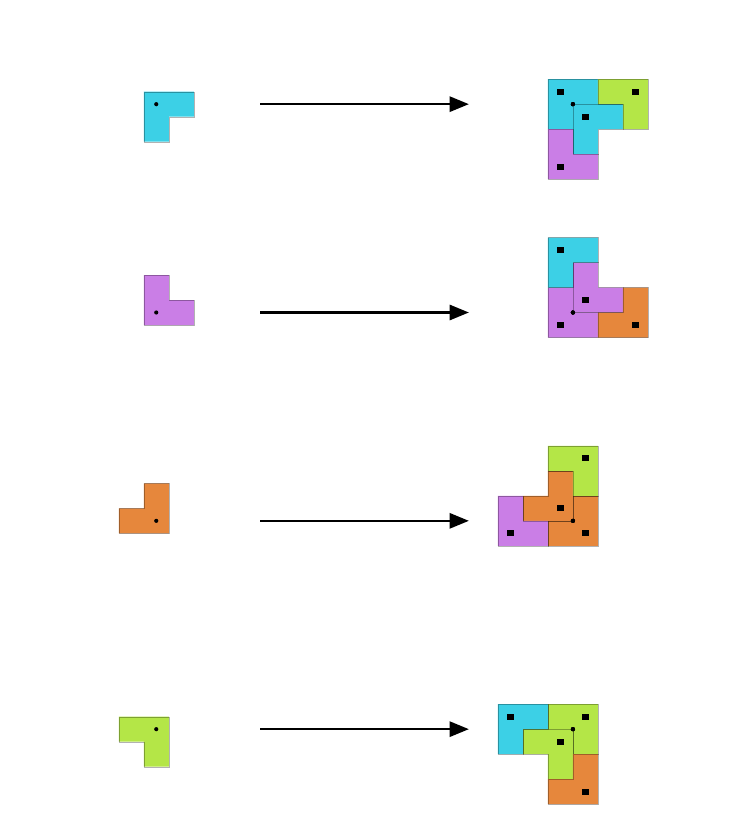}
        \caption{On the left, the prototiles \(\tau_1, \tau_2, \tau_3, \tau_4\)
        with their origins marked.
        On the right,
        the sets \(\{T_x \tau_k \mid k \in [4], x \in \Delta_{kj}\}\),
        with the origin of each tile marked as a square,
        and the origin of \(D_2 \tau_j\) marked as a small circle.
        From the pictures, we see that \(\varrho(\bm{o}_j)\)
tiles \(D_2 \tau_j\) for \(j \in [4]\).
        The \(C_4\)-rotational symmetry is also visible.
    }
    \label{fig:chair-substitution-rule}
\end{figure}

\begin{figure}[h]
    \centering
    \includegraphics[width=0.8\textwidth]{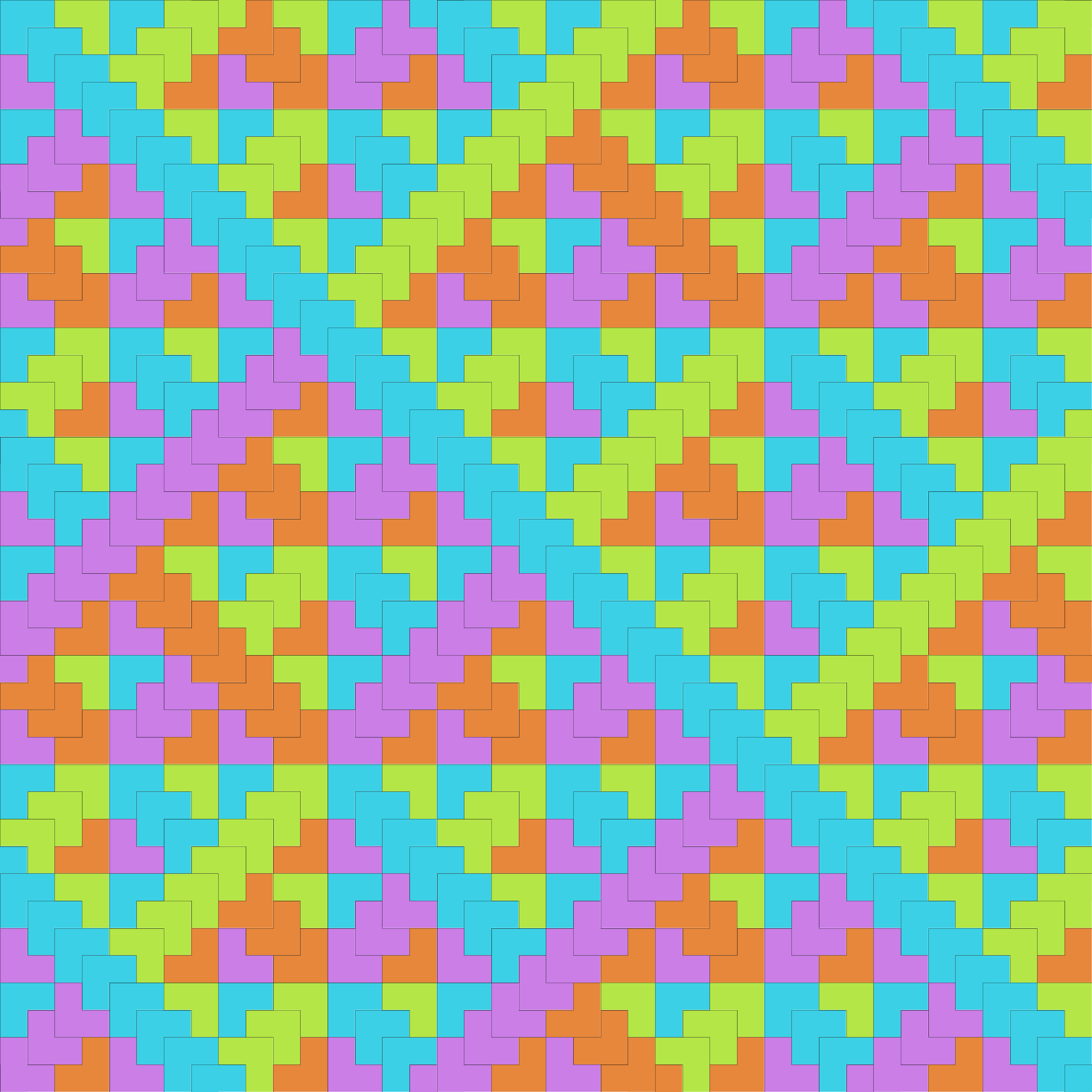}
    \caption{A chair tiling, i.e. a sample of the set \(\{T_x \tau_j \mid j \in [4], x \in \Phi_j\}\),
    where \(\bm{\Phi} = (\Phi_j)_{j \in [4]}\) is the vector point process associated to the chair substitution rule.}
    \label{fig:chair-tiling}
\end{figure}

\begin{example}\label{chair-rule-def}
We define the \emph{chair substitution rule} \(\cS\) on \(\bR^2\)
as the \(C_4\)-symmetric rule
given by Figure \ref{fig:chair-substitution-rule}.
We check that this is a primitive, stone substitution rule.

\begin{itemize}
    \item From the picture,
it is clear that \(\varrho(\bm{o}_j)\)
tiles \(D_2 \tau_j\) for all \(j \in [4]\)
with prototiles \(\tau_1, \tau_2, \tau_3, \tau_4\),
therefore \(\cS\) is stone with scaling constant \(\lambda = 2\)
and the prototiles \(\tau_1, \tau_2, \tau_3, \tau_4\)
are its canonical prototiles.
\item 
After applying the substitution rule twice,
all four orientations of the original prototile appear,
hence \(\cS\) is primitive.
\end{itemize}
Figure \ref{fig:chair-tiling} shows a chair tiling,
i.e. a tiling defined by an element of the associated substitution space \(\Omega_\cS\).

\end{example}

There is only one prototile up to rotation.
Therefore,
the spherical substitution matrix is
\(\Msph = (4)\).

\begin{corollary}
The vector point process associated to the chair rule is hyperuniform for constant weights.
\end{corollary}

The full substitution matrix
of the chair substitution rule
is given by \[
    \Mfull = \begin{pmatrix}
        2 & 1 & 0 & 1 \\
        1 & 2 & 1 & 0 \\
        0 & 1 & 2 & 1 \\
        1 & 0 & 1 & 2
    \end{pmatrix},
\]

as one can read from \Cref{fig:chair-substitution-rule}.
It has eigenvalues \(4,2,2,0\), with eigenvectors
\[
    v_1 = \begin{pmatrix} 1 \\ 1 \\ 1 \\ 1 \end{pmatrix}, \quad
    v_2 = \begin{pmatrix} 1 \\ i \\ -1 \\ -i \end{pmatrix}, \quad
    v_3 = \begin{pmatrix} 1 \\ -i \\ -1 \\ i \end{pmatrix}, \quad
    v_4 = \begin{pmatrix} 1 \\ -1 \\ 1 \\ -1 \end{pmatrix}.
\]

Then \Cref{thm:main-theorem}
gives us the following bounds.
\begin{corollary}
\label{corollary:chair-uniformity}
Let \(\cS\) be the chair substitution rule
and \(w \in \bC^4\).
Then \begin{enumerate}
    \item If \(w \in \linspan \{v_1, v_4\}\),
    then \(\hat{\eta}_w(B_r) \lessapprox r^4\) as \(r \to 0\).
    In particular, the point process associated to \(\cS\) is hyperuniform of Class I with these weights.
    \item If \(w \notin \linspan \{v_1, v_4\}\),
    then \(\hat{\eta}_w(B_r) \lessapprox r^2\) as \(r \to 0\).
\end{enumerate}
\end{corollary}
\begin{proof}
If \(w \in \linspan \{v_1, v_4\}\),
\(w\) is orthogonal to the eigenspaces of the eigenvalue \(2\),
therefore \(\beta_\perp(w) = 3 - \log_2 2 = 2\)
and \(\beta_\parallel(w) = 2 - \log_2 0 = \infty\).
Then \(\beta(w) = 2\),
so \Cref{thm:main-theorem} implies
\(\hat{\eta}_w(B_r) \lessapprox r^4\).

Otherwise, we have \(\beta(w) = \beta_\parallel(w) = 2 - \log_2 2 = 1\),
so \(\hat{\eta}_w(B_r) \lessapprox r^2\).
\end{proof}

Part \textit{(ii)}
does not prove hyperuniformity for these weights:
if the bound were sharp,
this would mean that the point process associated to the chair substitution rule
is not hyperuniform for these weights.

In the literature,
one often considers a variant of the chair substitution rule,
the \emph{block substitution rule} \citep*{limit-periodic},
which has the same substitution matrix
but where the prototiles are all squares:
as the substitution matrix is the same,
the bounds from \Cref{corollary:chair-uniformity} still hold.

For the block substitution rule, self-similar sets
are \emph{limit-periodic},
and \citet*{limit-periodic} computed the diffraction measure explicitly.
In this case,
we can prove that there are choices of weights
for which the point process associated to the block substitution rule is not hyperuniform.
This suggests that the same is true for the chair substitution rule.

\begin{theorem}
Let \(\cS_{block}\)
be the block substitution rule
and \(w = (1, 0, 0, 0)^\top\).
Then \(\hat{\eta}_w(B_r) \gtrsim r^2\) as \(r \to 0\).
\end{theorem}
\begin{proof}
By \Cref{lemma:uniformity-from-density},
it suffices to find a subset \(A \subset B_r\)
which is not a null set
and such that \(\hat{\eta}_w(D_{2^{-n}}A) \gtrsim 4^{-n}\):
as \(\cS_{block}\) is pure point diffractive,
it suffices to do this for a single point.

In fact, \citet*[Equation (25)]{limit-periodic} showed that
\[
    \hat{\eta}_w(2^{-n}(1,0)^\top) = \left| \frac{1}{4^n} \frac{1}{1 - \epsilon_n} \right|^2 \asymp 4^{-n}
\]
as \(n \to \infty\),
where \(\epsilon_n = e^{-2\pi i 2^{-n}}\).
This completes the proof.
\end{proof}

\subsubsection{Other examples}
The \emph{Tilings Encyclopedia} \citep*{encyclopedia}
is a large compendium of interesting tilings
including many coming from symmetric substitution rules.
As long as they are primitive and stone,
one can apply Theorem \ref{thm:main-theorem}
to try to prove hyperuniformity for many of these tilings:
in order to do this, we only need to find their spherical matrix.
All the information we need can be obtained from the pictures in the encyclopedia.

\begin{theorem}\label{thm:example-list}
The vector point processes associated to the following substitution rules are hyperuniform for any symmetric choice of weights:
\begin{enumerate}
    \item Penrose \citep*{thesis} (see \Cref{fig:example-thm-tilings-intro}).
    \item \href{https://tilings.math.uni-bielefeld.de/substitution/danzers-7-fold/}{Danzer's 7-fold} \citep*{encyclopedia_danzers_7fold} (see \Cref{fig:example-thm-tilings-intro}).
    \item \href{https://tilings.math.uni-bielefeld.de/substitution/binary}{Godr\`eche--Lan\c{c}on--Billard} \citep*{encyclopedia_binary} (modified as in \citep*[Fig.~3]{Baake_2019_renormalisation}, see \Cref{fig:example-thm-tilings-intro}).
    \item \href{https://tilings.math.uni-bielefeld.de/substitution/cap/}{CAP} \citep*{encyclopedia_cap}.
    \item \href{https://tilings.math.uni-bielefeld.de/substitution/watanabe-ito-soma-12-fold/}{Watanabe--Ito--Soma 12-fold} \citep*{encyclopedia_watanabe_ito_soma_12fold}.
    \item \href{https://tilings.math.uni-bielefeld.de/substitution/ammann-a3/}{Ammann A3} \citep*{encyclopedia_ammann_a3}.
\end{enumerate}
\end{theorem}
\begin{proof}
See the given references for the full definitions of these substitution rules.
Note that, for rules (iii)--(vi),
the prototiles depicted in the encyclopedia are not the canonical prototiles.
However, we still know the rules are stone,
because they define self-similar Delone subsets.
For the Godr\`eche--Lan\c{c}on--Billard,
the canonical prototiles have fractal boundary
and have been computed by \citet*{glb}.
One can also check that they are primitive from the full substitution matrices.

In order to prove the associated vector point processes are hyperuniform,
we need to compute the spherical substitution matrices of the substitution rules
and compare their largest non-PF eigenvalue \(\mu_2\)
with the scaling constant \(\lambda\).
We do this in \Cref{tab:examples-summary}.
If the scaling constant is not written in the encyclopedia,
we can compute it as \(\lambda = \sqrt{\lambda_{PF}}\),
where \(\lambda_{PF}\) is the PF eigenvalue.
For all of the examples, we see that \(\abs{\mu_2} < \lambda\),
which means they are hyperuniform for any choice of symmetric weights
by \Cref{corollary:symmetric-substitution-systems-uniformity}. \qedhere
\end{proof}

\renewcommand{\arraystretch}{1.5}
\begin{table}[ht]
\centering
\begin{tabular}{l c c c}
\hline
Rule & $\lambda$ & $\Msph$ & $\mu_2$ \\
\hline
Penrose & $\tau$ & $\begin{psmallmatrix}2&1\\1&1\end{psmallmatrix}$ & $-\tau^{-2}$ \\
\href{https://tilings.math.uni-bielefeld.de/substitution/danzers-7-fold/}{Danzer's 7-fold} & $\sqrt{1 + \frac{\sin(\frac{2\pi}{7})}{\sin(\frac\pi{7})}}$ & $\begin{psmallmatrix}2&1&2\\1&5&3\\2&3&3\end{psmallmatrix}$ & ${\sim}2.1$ \\
\href{https://tilings.math.uni-bielefeld.de/substitution/binary}{Godr\`eche--Lan\c{c}on--Billard} (modified) & \(\tfrac12(5+\sqrt{5}) = {\sim}3.6\) & \(\begin{psmallmatrix}10 & 5 \\ 5 & 5\end{psmallmatrix}\) & \(\tfrac12(15-5\sqrt{5}) = {\sim} 1.9\) \\
\href{https://tilings.math.uni-bielefeld.de/substitution/cap/}{CAP} & \(\frac{3 + \sqrt{5}}{2}\) & \(\begin{psmallmatrix}
    3 & 1 & 1 & 2 \\ 1 & 0 & 0 & 0\\ 3 & 0 & 1 & 1 \\ 3 & 0 & 2 & 3
\end{psmallmatrix}\)  & 1 \\
\href{https://tilings.math.uni-bielefeld.de/substitution/watanabe-ito-soma-12-fold/}{Watanabe--Ito--Soma 12-fold} & \(2 + \sqrt{3}\) & \(\begin{psmallmatrix}
    7 & 8 & 16 \\ 2 & 3 & 6 \\ 2 & 2 & 4
\end{psmallmatrix}\) & \(7 - 4 \sqrt{3}\) \\
\href{https://tilings.math.uni-bielefeld.de/substitution/ammann-a3/}{Ammann A3} & \(\frac{3 + \sqrt{5}}{2}\)  & \(\begin{psmallmatrix}
    1 & 3 & 2 \\ 1 & 1 & 0 \\ 0 & 1 & 1
\end{psmallmatrix}\) & 1 \\
\hline
\end{tabular}
\caption{For each rule in \Cref{thm:example-list}, we list its scaling constant,
spherical substitution matrix, and largest non-PF eigenvalue.
In all of these cases, we have \(\abs{\mu_2} < \lambda\),
so the associated point process is hyperuniform for symmetric weights.
Note that both \(\lambda\) and \(\mu_2\) can be computed from the spherical substitution matrix alone.
}
\label{tab:examples-summary}
\end{table}
Furthermore, by analyzing the full substitution matrices,
we can obtain more precise bounds using \Cref{thm:main-theorem}:
in \citep*{thesis} this was done for the Penrose and Danzer rules.
\begin{theorem}
\leavevmode
\begin{enumerate}
    \item For the Penrose substitution rule and any symmetric \(w\),
    we have \(\hat{\eta}_w(B_r) \lesssim r^{\alpha}\)
    as \(r \to 0\),
    where \(\alpha = {\sim}2.85\)
    \item For Danzer's 7-fold substitution rule and any symmetric \(w\),
    we have \(\hat{\eta}_w(B_r) \lesssim r^{\alpha}\)
    as \(r \to 0\),
    where \(\alpha = {\sim}3.1\)
\end{enumerate}
\end{theorem}

 
\bibliography{lit}
\bibliographystyle{plainnat}
	

\end{document}